\documentclass[a4paper,twoside]{article}
\usepackage[latin1]{inputenc}
\usepackage[T1]{fontenc}
\usepackage{amssymb}
\usepackage[final]{lpretex}
\usepackage{title}
\usepackage[all,dvips]{xy}
\usepackage{amscd}
\EndOfDump
\usepackage{url}
\usepackage{htmlref}

\def\maxid{{\frak{m}}}

\newcommand\tH{\widetilde H}

\newcommand\depth{\operatorname{depth}}

\newcommand\tsE{\widetilde{\Cal E}}
\newcommand\tsF{\widetilde{\Cal F}}

\newcommand\oF{\overline{\mathbf{F}}}

\newcommand\cs{\mathbf{s}}

\def\W{\cansymb{W}}

\def\Gr{\operatorname{Gr}}

\newcommand\domU{\mathbf{U}}

\renewcommand\contract{\mathop{\mbox{\rule{5pt}{.6pt}\rule{.6pt}{6pt}}\,}}

\newcommand\tot{\mathbf{s}}

\DocVersion[Wed May 28 2003]{1}
\DocVersion[Mod DEl-Ill argument;Mon Feb 23 2004]{2}

\begin{document}
\title{On non-liftable Calabi-Yau threefolds}

\author{T. Ekedahl}
\address{Department of Mathematics\\
 Stockholm University\\
 SE-106 91  Stockholm\\
Sweden}
\email{teke@math.su.se}

\subjclass{Primary 14J32; Secondary 14B12, 14J28, 14F40}

\begin{abstract}
Only two ways to construct non-liftable Calabi-Yau threefolds are currently
known, one example by Hirokado (cf.\ \cite{hirokado99::ae+calab+yau}) and one
method of Schr\"oer (cf.\ \cite{schroeer03::some+calab+yau+witt}).
This article computes some cohomological invariants of these examples
of non-liftable Calabi-Yau threefolds, in particular it computes their
mini-versal deformations. One conclusion is that their mixed characteristic
mini-versal deformation spaces are actually smooth over the characteristic $p$
base field. Furthermore, a new family, constructed in the spirit of Schr\"oer,
is introduced and the same calculations are performed for it.
\end{abstract}
\maketitle

There are only a small number of known example of non-liftable Calabi-Yau
manifolds. There is a $3$-fold in characteristic $3$ constructed by Hirokado
(cf.\ \cite{hirokado99::ae+calab+yau}) and two examples of Schr\"oer (cf.\
\cite{schroeer03::some+calab+yau+witt}), one in characteristic $3$ and one (or
rather a family) in characteristic $2$. In view of the current scarcity of
examples it seems interesting to study their properties and particularly, in view
of the fact that Calabi-Yau varieties in characteristic $0$ are unobstructed,
their deformation theory. In this paper we shall indeed make a study of their
cohomological invariants as well as determining their miniversal
deformations. For Hirokado's example we obtain in particular the following result.
\begin{Theorem}[A]\label{Rigidity}
The Hirokado threefold $X$ is arithmetically rigid, i.e., for a $A \in C_{\Z_3}$ and
a deformation of $X$ over $A$ we have that $3A=0$ and the deformation is of the
form $X\times\Sp A$.
\end{Theorem}
\begin{remark}
This implies arithmetic rigidity over any (perfect) field $\k$ of characteristic $3$
as the hull commutes with base change $\W(\F_3) \to \W(\k)$. Alternatively, the
proof nowhere makes use of the fact that the base field is $\F_3$.
\end{remark}
The main observation used in the proof of this theorem is that $X$ has a
manifestation as a Deligne-Lusztig type variety associated to the
Grassmannian $\Gr(2,4)$. This presents $X$ as a small resolution of a complete
intersection of type $(2,4)$ in $\P^5$ which allows for the computation of all
its Hodge numbers.  This proves rigidity in characteristic $3$. To show
arithmetic rigidity we use an idea of Serre (cf. \cite{serre61::exemp}) and show
that if the variety lifted to $\W/9$ then also a representation of a group of
automorphisms of the variety would lift which it doesn't.

As for Schr\"oer's examples we shall show first that there is a third family of
the same type, namely admitting a smooth pencil of supersingular (in the sense
of having $\rho=22$) K3-surfaces. We then show the following result.
\begin{Theorem}[B]\label{Pencil}
A versal $\W(\k)$-deformation space for a Calabi-Yau threefold in positive
characteristic which is the total space of a smooth pencil of supersingular
K3-surfaces is a smooth $\k$-scheme.
\end{Theorem}
The proof consists first of a determination of the dimension of the tangent
space of a miniversal deformation combined with the fact that that dimension can
be realised by an actual deformation in characteristic $p$. This shows that a
versal $\k$-deformation is smooth. To then show that the $\W$-deformation space
is no larger we shall use the easily established fact that if the mod $p$
deformation space is smooth and the $\W$-deformation space is not killed by $p$
then there is a deformation over $\W/p$. This will almost contradict a theorem
of Deligne and Illusie as we shall see that the Hodge to de Rham spectral
sequence does not degenerate. The result of Deligne and Illusie is not
applicable as it stands in the case when we are dealing with threefolds in
characteristic two but we shall conclude using some additional reasoning. In
order to perform the required calculations as well as to show the existence of a
third family of examples we shall use Ogus' theory of periods of supersingular
K3-surfaces (cf.\ \cite{Og79}).

More precise results on the dimension of the miniversal deformation will be
given in section \ref{sec:pencils}. Just as in the case of Hirokado's example we
shall also compute the Hodge-cohomology and contrary to Hirokado's example we
compute the de Rham-cohomology.
\begin{section}{The cohomology of complete intersections}

In the next section we shall need an extension of some of the computations of
\cite{deligne73::cohom} to not necessarily smooth complete intersections. We
begin by a preliminary result on the wedge filtration of some not necessarily
locally free modules. Let us recall that if
\begin{displaymath}
\shex MNP
\end{displaymath}
is a short exact sequence of $\sO$-modules on a scheme $X$ then the
\Definition{wedge filtration} on $\Lambda^n N$ is given by the images of the
maps $\Lambda^iM\Tensor\Lambda^{n-i}N$ given by wedge multiplication.
\begin{lemma}\label{wedge filtration}
Suppose that $X$ is a noetherian scheme, $Z$ a closed subset and $\k$ an integer
such that $\depth_x\sO_{X,x}
\ge k$ for all $x \in Z$. Suppose that
\begin{displaymath}
\shex{\sE}{\sF}{\sM}
\end{displaymath}
is a short exact sequence with $\sE$ and $\sF$ locally free and $\sM$ locally
free outside of $Z$. Set $n:=\rk\sF-\rk\sE$.

\part[i] The successive quotients of the saturation of the wedge filtration of $\Lambda^m\sF$ are
isomorphic to $\Lambda^i\sM\Tensor\Lambda^{m-i}\sE$ for $m \le k$.

\part[ii] For every $x \in Z$ we have $\depth_x\Lambda^m\sM \ge k-m$.

\part[iii] If $k \ge 2$ there is a map $\Lambda^n\sM \to \det\sM$ that is an
isomorphism outside of $Z$, where $\det\sM$ is the determinant of any locally
free resolution of $\sM$ so that $\det\sM = \det\sF\Tensor(\det\sE)^{-1}$.

\part[iv] If $k -m \ge 2$ then wedge multiplication composed with the map $\Lambda^n\sM \to
\det\sM$ induces an isomorphism $\Lambda^{m}\sM \liso \sHom(\Lambda^{n-m}\sM,\det\sM)$.


\part[v]\label{Tor independence} Suppose that \map{\pi}XS and \map f{S'}S are maps of schemes
such that $\pi(x)=s=f(s')$. Assume that 
$\Tor^i_{\sO_{S,s}}(\sO_{S',s'},\sO_{X,x})=\Tor^i_{\sO_{S,s}}(\sO_{S',s'},\sM_{X,x})=0$
for all $i>0$. Then
$\Tor^i_{\sO_{S,s}}(\sO_{S',s'},\Lambda^m\sM_{X,x})=0$ for $m \le k$.
\begin{proof}
To begin with, note that, by
\cite[Prop.~III:3.3]{grothendieck68::cohom+lefsc+sga}, 
$\depth_x\sO_{X,x} \ge k$ for all $x \in Z$ if and only if $\sH_Z^i(\sO)=0$
for $i<k$ and that similarly $\depth_x\Lambda^m\sM\ge k-m$ if and only if
$\sH_Z^i(\Lambda^m\sM)=0$ for $i < k-m$. The case $m=0$ follows
from the assumptions. For $m>0$ we consider the wedge filtration $\{0\}\subseteq\sG^0\subseteq \sG^1 \subseteq
\dots\sG^m=\Lambda^m\sF$ and we have
\begin{itemize}
\item For $0 \le r \le m$ there is a surjective map
$\Lambda^{m-r}\sE\Tensor\Lambda^r\sM \to \sG^{r}/\sG^{r-1}$.

\item When $r=m$ this map is an isomorphism so that $\sG^m/\sG^{m-1}=\Lambda^m\sM$.
\end{itemize}
Note furthermore that
$\sH_Z^i(\Lambda^{m-r}\sE\Tensor\Lambda^r\sM)=\Lambda^{m-r}\sE\Tensor\sH_Z^i(\Lambda^r\sM)$
so that by induction $\sH_Z^i(\Lambda^{m-r}\sE\Tensor\Lambda^r\sM)=0$ if $r < m$
and $i \le k-m$. Outside of $Z$ $\sM$ is locally free and then it is well-known
that $\Lambda^{m-r}\sE\Tensor\Lambda^r\sM \to \sG^{r}/\sG^{r-1}$ is an
isomorphism. Hence its kernel has support in $Z$ and as
$\sH_Z^0(\Lambda^{m-r}\sE\Tensor\Lambda^r\sM)=0$ the kernel is actually
zero. Hence we have short exact sequences
\begin{displaymath}
\shex{\sG^{r-1}}{\sG^{r}}{\Lambda^{m-r}\sE\Tensor\Lambda^r\sM}
\end{displaymath}
which imply that $\sH_Z^i(\sG^{m-1})=0$ when $i \le k-m+1$ and in any case
$\sH_Z^i(\sG^{m})=0$ when $i < k$ as $\sG^m=\Lambda^m\sF$ is locally free. Then
using the short exact sequences
\begin{displaymath}
\shex{\sG^{m-1}}{\sG^{m}}{\Lambda^m\sM}
\end{displaymath}
we conclude that $\sH_Z^i(\sG^{m-1})=0$ when $i \le k-m$.

Outside of $Z$ we have an isomorphism $\Lambda^n\sM = \det\sM$ and by the depth
conditions on $\sO$, $\det\sM=j_*\det\sM_U$, where \map j{X\setminus Z}X, so that
the map of \DHrefpart{iii} is simply the adjunction map. As for \DHrefpart{iv},
again by the depth condition $\sHom(\sN,\sP)=j_*\sHom(\sN_U,\sP_U)$ for any
$\sO$-module $\sM$ and any locally free $\sO$-module $\sP$. This reduces
\DHrefpart{iv} to the same statement over $U$ where it is clear.


Finally, \DHrefpart{v} follows easily from \DHrefpart{i} and induction.
\end{proof}
\end{lemma}
We shall now prove an extension of \cite[Prop.~1.3]{deligne73::cohom} which
will allow singularities. As we shall not need it, we will however dispense with the
auxiliary sheaf $\sF$ of [loc.~cit.] (the proof will not go through without some
conditions on $\sF$).
\begin{proposition}\label{cohomology vanishing}
Let \map{f}XS be a smooth and proper morphism, $\sW$ a locally free sheaf on $X$
of rang $c$, $d$ an integer $>0$ and \map s{\sW}{\sO_X} a section of
$\sW^*$. Suppose that $S$ is noetherian and that
\begin{enumerate}
\item The singular locus of the map $H \to S$ has fibrewise codimension $\ge k$,
where $H$ is the subscheme of $X$ defined by the vanishing of $s$.

\item If the $s_i$ are the coordinates for $s$ in a local basis of $\sW\check{}$
then the $s_i$ is a regular $\sO_X$-sequence.

\item For any sequence $k_i$ of integers, not all of which are zero, one has
\begin{displaymath}
R^if_*(\Tensor_i\Lambda^{k_i}\sW\Tensor\Omega_{X/S}^j)=0
\end{displaymath}
for $i+j < d$ and $j < k$.
\end{enumerate}
Then
\begin{enumerate}
\item $R^if_*\Omega_{X/S}^j\riso R^if_*\Omega_{H/S}$ for $i+j < d-c$ and
$j < k$.

\item $R^if_*\Omega_{X/S}^j\hookrightarrow R^if_*\Omega_{H/S}$ for $i+j < d-c$ and
$j < k$.
\end{enumerate}
\begin{proof}
Just as in \cite[Prop.~1.3]{deligne73::cohom} we will use the exact sequence
\begin{displaymath}
\shex{\sW_H}{\Omega_{X/S}^1}{\Omega_{H/S}^1}.
\end{displaymath}
As $H$ is locally defined in $X$ by a regular sequence $H$ is flat over $S$ and
as the above exact sequence remains exact over fibres of $f$ so is
$\Omega_{H/S}^1$. Hence the $\Tor$-conditions of (\ref{wedge filtration}) are
fulfilled and the formation of $\Omega_{H/S}^m$ commutes with base change for $m
\le k$. Using this, localising and using Grothendieck's theorem on formal
functions we reduce to the case when $S$ is Artinian and in particular
Cohen-Macaulay. This means that the depth conditions of Lemma \ref{wedge
filtration} are fulfilled. Using it, the proof of
\cite[Prop.~1.3]{deligne73::cohom} goes through essentially unchanged.
\end{proof}
\end{proposition}
Just as in \cite[Thm.\ 1.5]{deligne73::cohom} we may combine the proposition with
duality to obtain results on the cohomology above the diagonal. However, because
of the presence of singularities, the results one would obtain would be much
weaker and in fact too weak to be of any use to us. We therefore restrict
ourselves to giving the immediate application of what we have obtained so far to
complete intersections in projective space.
\begin{corollary}\label{hodge numbers}
Let $Y$ be an $n$-dimensional complete intersection in $\P^r$ with only isolated
singularities. For the Hodge numbers $h^{ij}_Y:=\dim H^j(Y,\Omega_Y^i)$ we have
that $h^{ij}=\delta_{ij}$ when $i+j<n$. Furthermore, the Euler characteristic of
$\Omega^i_Y$ for $i<n$ is the same as that for a non-singular complete
intersection (cf.\ \cite[Thm.\ 2.3]{deligne73::cohom}).
\begin{proof}
This is a direct consequence of the proposition and for the statement on the
Euler characteristic by its proof.
\end{proof}
\end{corollary}
\end{section}
\begin{section}{Hirokado's variety}

We shall now consider the deformation theory of Hirokado's
(\cite{hirokado99::ae+calab+yau}) example of a non-liftable Calabi-Yau threefold
in characteristic $3$. Before that we will need some preliminaries on
foliations. Note that we shall find it convenient in this section to define
projective space as consisting of lines rather than hyperplanes of a vector
space.

Let \map fXY be a purely inseparable finite flat map of degree $p$ between
smooth varieties over a perfect field $\k$ of characteristic $p$ and let $\Cal
E$ be the kernel of the differential of $f$. Suppose that \map s{S := \Sp
\k[\delta]}Y is a tangent vector of $Y$ and consider the fibre product of it and
$f$. It gives a finite flat map $T \to S$ of degree $p$ and as $f$ locally is
given by the extraction of a $p$'th root of a function $T$ is of the form
$R=\k[\delta][t]/(t^p-s)$ with $s \in \k[\delta]$. Up to addition of a $p$'th
root to and scaling of $s$ we have only two possibilities: Either $s = \delta$
or $s=0$. In the latter case $T \to S$ has a section and the tangent vector $s$
lifts to a tangent vector of $X$. In the former case it doesn't and $R$ is
isomorphic to $\k[t]/(t^{2p})$ with $\delta$ mapping to $t^p$.

Suppose conversely that we have a map \map t{T := \Sp\k[t]/(t^{2p})}X and we are
looking for conditions that the composite $T \to X \to Y$ factor through $T \to
S$, where the latter map is given by $\delta \to t^p$. This happens precisely
when the equivalence relation on $T$ induced by $T \to S$ maps into the
equivalence relation on $X$ induced by $X \to Y$. The latter equivalence
relation is induced by the foliation $\Cal E$ so the condition translates into
requiring that for any local section of $\Cal E$ around the image point of the
image of the closed point of $T$ there should be a relative vector field of $T
\to S$ that is compatible with that section. As $\Cal E$ is locally generated by
one element it is enough to verify that on such a generator. In particular the
image of $T$ should be tangent to the foliation but as we shall see this is not
enough in general.


Let us now momentarily consider the general situation of a smooth variety (in
any characteristic) and a sub-bundle $\Cal E$ of rank $r$ of the tangent bundle
that is closed under commutators. Being a sub-bundle it gives a section $e$ of
the Grassmannian bundle $\Cal G$ of $r$-dimensional subspaces of the tangent
bundle. Being, as we shall see, somewhat analogous to the Gauss map we will call
it the \Definition{Gauss section} associated to $\Cal E$.
\begin{lemma}
Let $\Cal E$ be a rank $r$ subsheaf of the tangent sheaf of a smooth variety $V$
that is closed under commutators, let $\Cal G$ be the fibration of
$r$-dimensional subspaces of the tangent sheaf and $e$ the Gaussian section of
$\Cal G$ associated to $\Cal E$. Suppose we have chosen coordinate functions
$x_1,\dots,x_n$ on $U$ and use them to trivialise the tangent bundle and hence
$\Cal G$. Then the differential of the composite $f$ of $e$ and the projection
$\Cal G \to \Gr(r,n)$ vanishes on $\Cal E$.
\begin{proof}
Pick a geometric point $s$ of $V$. After a linear change of coordinates we may
assume that $\Cal E_s$ is spanned by $\partial/\partial
x_1,\dots,\partial/\partial x_r$. For a tangent vector $v \in T_{V,s}$ we claim
that we may compute its image under the differential of $f$ as follows. Lift $v$
to a vector field $D$ on $U$. Take any vector field $E \in \Cal E$ and consider
the commutator $[D,E]$. Evaluate it at $s$ and consider the coefficients in
front of $\partial/\partial x_{r+1},\dots,\partial/\partial x_n$. This gives a
linear map $\Cal E_s \to T_{V,s}/\Cal E_s$ which is the differential of
$v$. That this is so is easily seen as the value at $s$ of any $E \in \Cal E$
has zero coefficients in front of $\partial/\partial
x_{r+1},\dots,\partial/\partial x_n$. The lemma now follows directly as $\Cal E$
is closed under commutators.
\end{proof}
\end{lemma}
Let us now consider the special case when we have an open subset $V$ of
$n$-dimensional projective space and a sub-bundle $\sE$ of $T_V$ closed under
commutators. We may consider the inverse image $\tsE$ of $\sE$ in $T_{\tilde
V}$, where $\tilde V$ is the inverse image of $V$ under the natural map
$\A^{n+1}\setminus\{0\}$. Then $\tsE$ will also be closed under commutators and
hence has a Gauss section. However, $T_{\tilde V}$ has a trivialisation coming
from linear vector fields and hence we may compose the Gauss section with
projection of the Grassmann bundle onto $\Gr(r+1,n+1)$ giving a map
\map{f}{\tilde V}{\Gr(r+1,n+1)}. This map is invariant under the $\mul$-action
on $\tilde V$ and so descends to a map $V \to \Gr(r+1,n+1)$. Working only on $V$
it may be described as follows: Take the inverse image of $\sE$ under the
natural surjection $\Cal O(1)^{n+1} \to T_V$ and consider the section thus
obtained of the Grassmann bundle of $\Cal O(1)^{n+1}$. This bundle is trivial so
we may compose with the projection to $\Gr(r+1,n+1)$.
\begin{proposition}
With the above notations the kernel of the differential of $f$ contains $\Cal
E$. In particular, if $\Cal E$ is integrable and the quotient by it exists (as an
algebraic variety) then $f$ factors through this quotient.
\begin{proof}
This follows immediately from the lemma applied to $\tsE$.
\end{proof}
\end{proposition}
We shall call this map the \Definition{Gauss map} of  $\Cal E$.

Let us now put ourselves in characteristic $p$ and consider the homogeneous
vector field $D := \sum_{i=0}^nx_i^p\partial/\partial x_i$, where we throughout
this section will assume that $n \ge 3$ ($n=1$ is trivial while $n=2$ can be
analysed in a similar fashion). It and the Euler vector field
$\sum_{i=0}^nx_i\partial/\partial x_i$ generates a ($p$-)integrable foliation on
$\A^{n+1}$ (which is singular at the $\F_p$-rational points) and descends to a
foliation $\Cal E$ on $\P^n$ which is easily seen to be rationally generated by
the vector field $\sum_{i=1}^n(x_i^p-x_i)\partial/\partial x_i$ considered by
Hirokado in (\cite{hirokado99::ae+calab+yau}). It has singularities at
$\P^n(\F_p)$. We shall study the Gauss map of this foliation (and at the same
time reprove some of Hirokado's results).
\begin{lemma}\label{Hirokado vector field}
\part $\Cal E$ is invariant under $\PGL_{n+1}(\F_p)$.

\part[ii] The Gauss map $\P^n\setminus \P^n(\F_p) \to \Gr(2,n+1)$ associated to
$\Cal E$ extends to a regular map from the blowing up $\tilde\P^n$ of $\P^n$ at
$\P^n(\F_p)$. Restricted to an exceptional divisor this extension is an
embedding with image the lines that contain the blown up point. In
particular the foliation given by $\Cal E$ is smooth on all of $\tilde\P^n$.

\part[iii] The extended Gauss map $\tilde\P^n \to \Gr(2,n+1)$ contracts the (strict
transforms of) $\F_p$-rational lines in $\P^n$.

\part[iv] Two points of $\tilde\P^n$ map to the same point of $\Gr(2,n+1)$ if and
only they lie on the strict transform of the same $\F_p$-rational line.

\part[v] Let $X$ be the quotient of $\tilde\P^n$ by $\Cal E$. Then the induced map
$X \to \Gr(2,n+1)$ is an immersion outside of the images of the strict transform
of the $\F_p$-rational lines.
\begin{proof}
The first part is clear as the Euler vector field is invariant under all linear
transformations and $\sum_{i=0}^nx_i^p\partial/\partial x_i$ transforms like the
Euler vector field for transformations in $\GL_{n+1}(\F_p)$.

For the second part we notice that in homogeneous coordinates the Gauss map is
given by associating to $(x_0\co\dots\co x_n)$ the space spanned by the $2\times
(n+1)$-matrix
\begin{displaymath}
\left(
\begin{array}{cccc}
x_0^p&x_1^p&\dots&x_n^p\\
x_0  &x_1  &\dots&x_n
\end{array}
\right)
\end{displaymath}
if it is $2$-dimensional. By the first part we need only look at the point
$(1\co0\co\dots\co0)$ and furthermore we still have enough of transformations
left so that it is enough to look at only one of the patches of the blowing
up. We may subtract the top row from the bottom, substitute $x_i = s_is_n$ for
$1 \le i < n$ and $s_n=x_n$, divide the bottom row by $s_n$ and finally put
$s_n$ equal to $0$. This shows that the restriction of the Plücker embedding to
the chosen patch maps $(s_1,s_2,\dots,s_{n-1})$ to the $2$-dimensional space
spanned by the vectors $(1,0,\dots,0)$ and $(0,s_1,\dots,s_{n-1},1)$. These are
clearly independent so that the map has no base points and equally clearly is an
embedding. Furthermore, we see that we get in this way exactly the lines that
pass through $(1\co0\co\dots\co0)$.

To continue with \DHrefpart{iii} let to begin with $x$ and $y$ be two distinct
non-$\F_p$-points in $\P^n$ and assume that the map to the same point in
$\Gr(2,n+1)$. This is equivalent to $(x,x^p)$ and $(y,y^p)$ spanning the same
$2$-dimensional vector space $V$. By assumption $x$ and $y$ are linearly
independent so they span this common space and hence $x^p$ and $y^p$ are
contained in $V$. This means that $V$ is stable under the $p$'th power and hence
is defined over $\F_p$ so that $x$ and $y$ lie on the same $\F_p$-rational
line. Conversely of $x$ and $y$ lie in the same $\F_p$-rational vector space
then so do $x^p$ and $y^p$ and hence $(x,x^p)$ and $(y,y^p)$ span the same
vector space.

Turning to \DHrefpart{iv} let now $x$ be a non-rational point $y$ a point on an
exceptional divisor and assume that they map to the same point. Using the action
of $\PGL_{n+1}(\F_p)$ we may assume that $y$ lies over $(1\co0\co\dots\co0)$ and
in the same patch as in the proof of \DHrefpart{ii} and so it's image is the
space spanned by $(1,0,\dots,0)$ and $(0,s_1,\dots,s_{n-1},1)$ for some
$s_i$. On the other hand $x$ maps to the space spanned by $x$ and $x^p$. Hence
by assumption $y=(x_1,\dots,x_n)$ and $y^p$ span a $1$-dimensional space and so
it's a multiple of an $\F_p$-rational point and we may assume it is
$\F_p$-rational. This forces the $s_i$ to be rational and hence the line that is
the common image is $\F_p$-rational and $x$ lies on a rational line and $y$ is
easily seen to be the intersection of the strict transform of that line and the
exceptional divisor.

As for \DHrefpart{v} let us start by considering a point $y$ that lies under a
geometric point $x$ of $\P^2$ which does not lie on a $\F_p$-rational line. We
identify tangent vectors to $\Gr(2,n+1)$ at the space $V$ spanned by $(x,x^p)$
first with linear maps $V \to \k^{n+1}/V$ and then with pairs of vectors $(u,v)$
of $\k^{n+1}/V$, where $u$ is the image of $x$ and $v$ of $x^p$. In this form a
tangent vector $u$ of $T_x\P^n$, which itself is a vector in $\k^{n+1}/x$, is
mapped to the pair $(u,0)$. This gives the image of the subspace of $T_yX$ of
vectors that lift to $\tilde\P^n$. Consider now the $\P^n$-point
$x+tx^p+t^{p+1}x^{p^2}$ of $S := \k[t]/(t^{2p})$. We confuse $x$ with a
generator of that line in $\k^{n+1}$ and hence may lift the situation to the
affine space and consider $\sum_ix_i^p\partial/\partial x_i$ as a local
generator of the foliation. This generator takes $x+tx^p+t^{p+1}x^{p^2}$ to
$x^{p}+t^px^{p^2}$ which on the other hand can be obtained by applying the
vector field $\partial/\partial t$ to $x+tx^p+t^{p+1}x^{p^2}$. Hence our descent
condition on this point is fulfilled and it descends to a tangent vector on
$X$. (Note that we are forced to add the term $t^{p+1}$ to make this work hence
we have given an example of the fact that it is not enough to have the point
being tangent to the foliation.) Let us now consider the image of the point
under the Gauss map. It is mapped to the space spanned by $x^p+t^px^{p^2}$ and
$x+tx^p+t^{p+1}x^{p^2}$ which is the same as the space spanned by
$x^p+t^px^{p^2}$ and $x$ and hence the tangent vector on $Y$ is mapped to
$(0,x^{p^2})$. As $x$ does not lie on a $\F_p$-rational line, $x^{p^2}$ does not
lie in $V$ and hence the image of the tangent space $T_yY$ in the tangent space
of the Grassmannian is $n$-dimensional. What remains to consider is a point on
an exceptional divisor away from an intersection with the strict transform of an
$\F_p$-rational line which is simpler than the one just considered and left to
the reader.
\end{proof}
\end{lemma}
We have thus shown that the Hirokado variety $X$ is a (very) small resolution of
singularities of its image in $\Gr(2,n+1)$. We shall now identify this image
with a Deligne-Lusztig style variety. Hence let $\sF$ be the closed subscheme of
$\Gr(2,n+1)$ of $2$-dimensional subspaces $V$ of $\k^{n+1}$ such that $V \cap
F^*V \ne \{0\}$, where $F$ is the relative Frobenius map and let $\tsF$ be the
scheme of partial flags $L \subset V \subset \k^{n+1}$ with $\dim L=1$ and $\dim
V=2$ such that $L \subset F^*V$. We have a forgetful map $\tsF \to \sF$ which is
an isomorphism outside of the $\F_p$-rational subspaces $V$.
The Gauss map $\P^n\setminus \P^n(\F_p) \to
\Gr(2,n+1)$ takes $x$ to the space $V$ spanned by $x$ and $x^p$. Then $x^p
\subseteq V\cap F^*V$ and hence the image of the Gauss map lies in $\sF$. On
the other hand, if we put $L:=x^p$ we have a flag $L \subset V$ with $L \subset
F^*V$ and hence we get a map $\P^n\setminus \P^n(\F_p) \to \tsF$ which is easily
seen to extend to $\tilde\P^n \to \tsF$. Birationally this is just the Gauss map
so it factors to give a map $X \to \tsF$.
\begin{proposition}
\part The map $X \to \tsF$ is an isomorphism.

\part $\sF$ is a Cohen-Macaulay variety and $\tsF \to \sF$ is a resolution of
singularities of $\sF$.

\part \label{normal bundle} The positive dimensional fibres of $\tsF \to \sF$
are $\P^1$'s with conormal bundle isomorphic to $\sO(1)^{n-1}$.
\begin{proof}
We consider one further scheme $\tsF'$ namely that of a flag $L' \subset V
\subset \k^{n+1}$ with $\dim V=2$ and $\dim L'=1$ and $F^*L' \subset V$. We map
$\tsF'$ to $\tsF$ by mapping $(L' \subset V)$ to $(F^*L' \subset V)$. We get a map
$\tsF' \to \P^n$ by forgetting $V$. It is clear that this map is an isomorphism
at the open subscheme where $L' \ne F^*L'$ as then $V=L'+F^*L'$ and that the
exceptional set maps to $\P^n(\F_p)$. It is equally clear that this map is just
the blowing up of the $\F_p$-rational points of $\P^n$ so that
$\tsF'=\tilde\P^n$. Going backwards it then also follows that the forgetful map
$\tsF' \to \sF$ is the same as the Gauss map.

It can be shown directly that $\tsF$ is smooth but it also follows from the first
part. It again follows from that part that $\dim\sF=\tsF=n$ and then the fact
that $\sF$ is Cohen-Macaulay follows from \cite[Thm.\ 14.3 (c)]{Fu84}.

Finally, it follows from (\ref{Hirokado vector field}), but can just as easily be
seen directly, that $\tsF \to \sF$ is an isomorphism outside of the points $V$ on
$\sF$ for which $F^*V=V$, i.e., the $\F_p$-rational points of $\Gr(2,n+1)$ and
that the map contracts exactly the curve of lines in $V$ for each such rational
$V$. To compute the conormal bundle we will use the modular interpretation of $\tsF$
and we therefore start with a short discussion on how to ``see'' the normal
bundle in that context.

Let thus $Y$ be a smooth variety over a field $\k$ and $Z \subset Y$ a smooth
subvariety. By construction maps $Z[\epsilon] \to Y$ that restricts to the
inclusion $Z \subset Y$ correspond to maps $\Omega^1_{Y|Z} \to \sO_Z$. We have
the conormal bundle sequence
\begin{displaymath}
\shex{\check{\sN}_{Y/Z}}{\Omega^1_{Y|Z}}{\Omega^1_Z},
\end{displaymath}
and the maps $\Omega^1_{Y|Z} \to \sO_Z$ that factor through $\Omega^1_{Y|Z} \to
\Omega^1_Z$ correspond to the composite of the trivial map $Z[\epsilon] \to Z
\subset Y$ and the automorphism of $Z[\epsilon] \to Z[\epsilon]$ corresponding
to the vector field $\Omega^1_Z \to \sO_Z$. Consequently let us consider the
presheaf $\sG$ that to each open subset $U \subseteq Z$ associates the set of
extensions of $Z \hookrightarrow Y$ to $Z[\epsilon]$ modulo the automorphisms of
$Z[\epsilon]$ that commute with multiplication by $\epsilon$ and induce the
identity on $Z$. Then the associated sheaf is the normal bundle of $Z$ in $Y$
(with the evident $\sO_Z$-module structure).

As all exceptional curves of $\tsF \to \sF$ are conjugate under the action of
$\GL_{n+1}(\F_p)$ it is enough to consider the rational plane spanned by
$e_1:=(1,0,\dots,0)$ and $e_2:=(0,1,\dots,0)$. We thus have $\P^1
\hookrightarrow \tsF$ given by $(x\co y) \mapsto (\Span{e_1,e_2},(x\co y))$. Let
$U$ be an open subset of $\P^1$ and consider an extension of $U \hookrightarrow
\tsF$ to $U[\epsilon]$. It is given by a line $L_\epsilon$ and a plane
$V_\epsilon$ with $L_\epsilon \subseteq V_\epsilon \cap F^*V_\epsilon$. Now,
$F^*V_\epsilon$ is the constant plane $V:=\Span{e_1,e_2}$ and thus in particular
$L_\epsilon \subset V$. By an $\epsilon$-automorphism of $U[\epsilon]$, which is
the identity on $U$ we can make $L_\epsilon$ constant, i.e., equal to $(x\co
y)$. This takes care of all automorphisms so that by assuming this we do not
need to further consider dividing out by automorphisms. The only variation left
is that of $V_\epsilon$. There is a unique basis of it of the form
$e'_1:=(1,0,\epsilon a_1,\dots,\epsilon a_{n-1})$ and $e'_2:=(1,0,\epsilon
b_1,\dots,\epsilon b_{n-1})$. It contains $L$ precisely when $xa_i+yb_i=0$ for
all $i$ or in other words that $(a_i,b_i)$ lies in the kernel of $(r,s) \mapsto
xr+ys$. The sheaf of such solutions is $\sO(-1)$ so that the normal bundle is
isomorphic to $\sO(-1)^{n-1}$.
\end{proof}
\end{proposition}
\begin{remark}
\part Let us recall the construction of Dickson of generators for the ring of
invariants for the action of $\GL_n(\F_p)$ on $\F_p[x_1,\dots,x_n]$. One first
considers the determinant
\begin{displaymath}
D_{n,i}(x_1,\dots,x_n) = \left| 
\begin{array}{llcl}
x_1^{p^{n+1}}&x_2^{p^{n+1}}&\cdots&x_n^{p^{n+1}}\\
x_1^{p^n}&x_2^{p^n}&\cdots&x_n^{p^n}\\
\multicolumn{4}{c}\dotfill\\
\widehat{x_1^{p^{i}}}&\widehat{x_2^{p^{i}}}&\cdots&\widehat{x_n^{p^{i}}}\\
\multicolumn{4}{c}\dotfill\\
x_1&x_2&\cdots&x_n
\end{array}
\right|
\end{displaymath}
and then a set of generators for the invariants consist of $d_i:=D_{n,i}/d_n$ for $1
\le i < n$ and $d_n^{p-1}$, where $d_n:=D_{n,n+1}$. The action of the vector
field $D$ on the $D_{n,i}$ is very simple, in fact all the rows of the matrix
consist of constants except for the first one and the first one is mapped into
the second row except when $i=1$. This gives $D(D_{n,i})=0$ if $i \ne 1$ and
$D(D_{n,1})=D_n^p$. This observation was the starting point for the
consideration of what is the Gauss map in the Hirokado case as its $2\times
2$-minors would be $D$-constants for the same reason. Note also that the action
of $D$ on the $d_i$ is equally pleasant; $D(d_i)=0$ if $i \ne 1$ and
$D(d_1)=d_n$ so that in particular the ring of $D$-constant and
$\GL_n(\F_p)$-invariant polynomials is the polynomial ring
$\F_p[d_1^p,d_2,\dots,d_n]$.

\part The derivation $D$ appears in another, from the geometric point of view
rather intriguing, context. Let $p$ be an odd prime. The subring of the
cohomology ring $H^*((\Z/p)^n,\F_p)$ generated by the elements $x_i$ of
$H^2((\Z/p)^n,\F_p)$ obtained by applying the Bockstein operator to the
projections $e_i$ of $H^1((\Z/p)^n,\F_p)=\Hom((\Z/p)^n,\F_p)$ is a polynomial
algebra $\F_p[x_1,\dots,x_n]$. The first Steenrod operation $P^1$ is a
derivation and takes $x_i$ to $x_i^p$ so equals $D$. The invariance of $D$ under
$\GL_n(\F_p)$ is then an immediate consequence of the fact that this group
operates by continuous maps on $K((\Z/p)^n,1)$. The situation for $p=2$ is
similar; the whole cohomology is a polynomial ring and $P^1$ is to be replaced
by $Sq^1$.
\end{remark}
\end{section}
\begin{section}{Hodge cohomology}

The last proposition of the previous section allows us to compute the higher
direct images of the sheaf of $i$-forms for all $i$ of the map $\tsF \to \sF$.
\begin{proposition}\label{resolution cohomology}
Let \map{\pi}{\tsF}{\sF} be the canonical map and let $E$ be the disjoint union
of the exceptional curves for $\pi$.

\part $R^i\pi_*\Omega^j_{\tsF}=0$  if $i > 1$ or $i = 1$ and $j \ne 1$ where $0
\le j \le n$.

\part The composite of the restriction map $R^1\pi_*\Omega^1_{\tsF} \to
R^1\pi_*\Omega^1_E$ and the trace map $R^1\pi_*\Omega^1_E \to \sO_S$, where $S$
is the (reduced) image of $E$ under $\pi$, is an isomorphism.

\part $\pi_*\Omega^i_{\sF}$ equals $j_*\Omega^i_U$, where \map jU{\sF} is the
inclusion of the non-singular locus.
\begin{proof}
It follows from (\ref{normal bundle}) that $\sI_E/\sI_E^2$ is isomorphic to
$\sO(1)^{n-1}$ and then $\sI_E^m/\sI_E^{m+1}$ is isomorphic to
$S^m(\sO(1)^{n-1})$. On the other hand we have the exact sequence
\begin{displaymath}
\shex{\sI_E/\sI_E^2}{\Omega^1_{\tsF}/\sI{\Omega^1_{\tsF}}}{\Omega^1_E}.
\end{displaymath}
From this, and the computation of cohomology of sums of line bundles on
$\P^1$ it follows immediately that
$H^i(\tsF,\sI^n\Omega^1_{\tsF}/\sI^{n+1}\Omega^1_{\tsF})=0$ for $i,n > 0$ and
that the map $\Omega^1_{\tsF}/\sI{\Omega^1_{\tsF}} \to \Omega^1_E$ induces an
isomorphism on $H^i$ for $i > 0$. The first two parts now follows from
Grothendieck's theorem on formal vs. ordinary direct images. The last statement
follows as the codimension of $E$ in $\tsF$ is $\ge 2$.
\end{proof}
\end{proposition}
In order to compute further we shall be forced to specialize to the case $n =
3$. We may embed $\Gr(2,4)$ in $\P^5$ by the Plücker embedding where it is the
quadric hypersurface given by the Plücker relation, the quadratic form
$q(x)$. If $b(x,y)$ is the bilinearisation of $q$ then the lines in $\P^3$ meet
if and only $b(a,b)=0$, where $a$ resp.\ are the Plücker coordinates of the two
lines. This means that $\sF$ is the complete intersection of
$\Gr(2,4)=\{q(x)=0\}$ and $\{b(x,x^p)=0\}$, hypersurfaces of degree $2$ and
$p+1$. We would now like to apply the previous proposition to compute the Hodge
numbers of $\tsF$ and using the fact that $\sF$ is a complete intersection. We
start this investigation with $\sF$.
\begin{remark}
The fact that $\sF$ is a complete intersection has an interesting consequence
when $p=3$. Starting with the Calabi-Yau threefold $X=\tsF$ we pass to $\sF$ by
contracting some curves, we can then deform $\sF$ to a smooth complete
intersection of type $(2,4)$ in $\P^5$ which again is a smooth Calabi-Yau
manifold. These steps are a well-known way to connect different Calabi-Yau
varieties and we see that Hirokado's threefold, which is not itself liftable,
can be connected to a liftable Calabi-Yau variety.
\end{remark}
\begin{subsection}{Linear algebra interlude}

We shall need some linear algebra preliminaries.

Let us recall (cf.\ \cite{roby63::lois}) that the free divided power algebra
$\Gamma(V)$ is homogeneous with degree $p$ component $\Gamma^p(V)$. It will
appear in our context as the dual of $S^p(V)$ is $\Gamma^p(V^*)$, indeed the
divided power algebra $\Gamma(V^*)$ is a Hopf algebra by mapping $\gamma_n(x)$
to $\sum_{i+j=n}\gamma_i(x)\tensor\gamma_j(x)$ which gives a Hopf algebra
structure on its (degreewise) dual $\Gamma(V^*)^*$. It equals $V$ in degree and
hence gives a map $S(V) \to \Gamma(V^*)^*$ which is the desired isomorphism.
This duality may also be described as identifying $S(V)$ with the algebra of
constant coefficient divided power differential operators on $\Gamma(V^*)$ or by
identifying $\Gamma(V^*)$ with constant coefficient differential operators on
$S(V)$. The maps between $\Gamma^*$'s that we will consider are of several
types. 
\begin{definition}\label{operators}
Let $V$ be a finite dimensional vector space over a field of characteristic $p
>0$.

\part For every $q \in S^n(V)$ the dual of multiplication
$S^m(V) \to S^{m+n}(V)$ by $q$ which gives a map $\Gamma^{m+n}(V^*) \to
\Gamma^m(V^*)$ which we shall call the \Definition{contraction operator}; it can also be
described as the action by $q$ considered as a divided power differential
operator. 

\part Picking out bihomogeneous components of the coproduct on $\Gamma(V)$ gives
maps $\Gamma^{i+j}(V) \to \Gamma^i(V)\Tensor\Gamma^j(V)$. We shall call them
\Definition{splitting operators}.

\part The maps $\Gamma^{pn}(V) \to F^*\Gamma^n(V)$, where $F$ is the
Frobenius map $\k \to \k$. They are dual to the $p$'th power maps $F^*S^n(V^*) \to
S^{pn}(V^*)$ and can described explicitly by the fact that they take
$\gamma_{pn}(v)$ to $\gamma_n(v)$. They will be called \Definition{Frobenius operators}.

\part In the special case when $V = \Lambda^2U$ we get a map $\Lambda^4U \to
\Gamma^2(\Lambda^2U)$; by picking out the $(1,1)$-component of the coproduct on
the exterior algebra gives a map $\Lambda^4U \to \Lambda^2U\Tensor\Lambda^2U$
and by the cocommutativity of the coproduct the image is contained in the
invariants under $\Sigma_2$, i.e., $\Gamma^2\Lambda^2U$. For any $k \ge 1$ we
then also get a map $\Gamma^k(\Lambda^4U) \to \Gamma^{2k}(\Lambda^2U)$ as the
composite of $\Gamma^k(-)$ applied to the map $\Lambda^4U \to
\Gamma^2(\Lambda^2U)$ composed with the natural map
$\Gamma^k(\Gamma^2(\Lambda^2U))$ that is induced by the map that to $u \in
\Gamma^2(\Lambda^2U)$ associates $\gamma_k(u) \in
\Gamma^{2k}(\Lambda^2U)$. These maps will be called the \Definition{copairing
maps}.
\end{definition}
Over a finite field we have a curious map from the rational points of projective
space on a vector space to some divided powers of it that will actually occur as a
map in a Leray spectral sequence.
\begin{definition-lemma}
Let $\F_q$ be finite field with $q$ elements and $V$ a finite-dimensional
$\F_q$-vector space. Let $\k$ be an integer divisible by $q-1$.

\part For a $1$-dimensional subspace $U \subseteq V$ and $0 \ne e \in U$,
$\gamma_k(e) \in \Gamma^kV$ is independent of $e$ and hence gives a well-defined
element that will be denoted $\gamma_k(U)$.

\part For an $n$-dimensional subspace $U \subseteq V$ we define $\gamma_k(U) \in
\Gamma^k\Lambda^nV$ by $\gamma_k(U) := \gamma_k(\Lambda^nU)$. We shall call it
the \Definition{$k$'th Pl\"ucker map}.
\begin{proof}
Another choice $e'$ has the form $\lambda e$ for $0 \ne \lambda \in \F_q$ and we
have $\gamma_k(e')=\lambda^k\gamma_k(e)$ and $\lambda^k=1$ as $q-1|k$.
\end{proof}
\end{definition-lemma}
\begin{proposition}\label{invariants} 
Let $V$ be four-dimensional and $q=p$, a prime, let $V_1
\subset V$ be a two-dimensional subspace and $G$ the group of linear
automorphisms of $V$ preserving $V_1$. 

\part The action of $G$ on
$\Gamma^{p-1}\Lambda^2V$ has exactly the space spanned by $\gamma_{p-1}(V_1)$ as
its space of invariants.  

\part\label{2-invariants} When $p=3$ the action on $\Gamma^{2(p-1)}\Lambda^2V$ has the space spanned
by $\gamma_{2(p-1)}(V_1)$ and the image of $\Gamma^{2(p-1)}(\Lambda^4V)$ under
the copairing map as its space of invariants.
\begin{proof}

We choose a complement $V_2$ to $V_1$ and start by considering the invariants
under the action of the subgroup of $G$ preserving also $V_2$. This subgroup is
isomorphic to $\GL_2(\F_p)\times\GL_2(\F_p)$ and we have
$\Lambda^2V=\Lambda^2V_1\Dsum V_1\Tensor V_2\Dsum\Lambda^2V_2$. This gives
\begin{displaymath}
\Gamma^k\Lambda^2V =
\Dsum_{i+j+m=k}\Gamma^i\Lambda^2V_1\Tensor\Gamma^j(V_1\Tensor V_2)\Tensor\Gamma^m\Lambda^2V_2.
\end{displaymath}


We have a cyclic subgroup of $\GL_2(\F_p)$ of order $p+1$ (corresponding to the
action of $\F_{p^2}^*$ on $\F_{p^2}$) and we may extend the coefficient field of
the representation such that we may find a basis $e_1,\dots,e_4$ of $V$ such
that $e_1,e_2 \in V_1$, $e_3,e_4 \in V_2$ and such that the $\F_{p^2}^*$ of the
first $\GL_2(\F_p)$-factor acts by $\lambda\cdot e_1=\lambda e_1$ and
$\lambda\cdot e_2=\lambda^{p}e_2$ and similarly for the second factor and $e_3$
and $e_4$. Put $e_{ij} := e_i\wedge e_j$. so that 
$\gamma_r(e_{13})\gamma_s(e_{14})\gamma_t(e_{23})\gamma_u(e_{24})$ for
$r+s+t+u=j$ form a basis for $\Gamma^j(V_1\Tensor V_2)$. Now $\lambda \in
\F_{p^2}^*$ in the first factor acts by multiplication of
$\lambda^{i(p+1)+r+s+(t+u)p}$ on $\Gamma^i\Lambda^2V_1\Tensor
\gamma_r(e_{13})\gamma_s(e_{14})\gamma_t(e_{23})\gamma_u(e_{24})\Tensor\Gamma^m\Lambda^2V_2$ while the same
$\lambda$ in the second factor acts by multiplication by
$\lambda^{m(p+1)+r+t+(s+u)p}$. Hence in order for it to be fixed by both factors
we must have
\begin{eqnarray*}
i(p+1)+r+s+(t+u)p &=& a(p^2-1)\\
m(p+1)+r+t+(s+u)p &=& b(p^2-1)
\end{eqnarray*}
for integers $a$ and $b$ and we also have $i+r+s+t+u+m=k$. Now, consider the case when $k=c(p-1)$,
$c=1,2$, so that $i+r+s+t+u+m=c(p-1)$. Adding the first two equations together
and using the last gives us $u-r=(a+b-c)(p+1)$. Now, we always have the
solutions $i=c(p-1),r=s=t=u=m=0$ and $m=c(p-1),i=r=s=t=u=0$ and for all others
we have $a,b > 0$. We also have $a(p^2-1)=i(p+1)+r+s+(t+u)p \le
(i+r+s+t+u+m)(p+1)=c(p^2-1)$ and similarly for $b$ so that $a,b \le c$ with
equality only if $r=s=t=t=u=0$ and $r=s=t=u=m=0$, i.e., the excluded cases so
that $a,b < c$. When $c=1$ this gives a contradiction.

When $c=1$ we thus have that the dimension of the
$\GL(2,\F_p)\times\GL(2,\F_p)$-invariants is at most two-dimensional. On the
other hand $\gamma_{p-1}(V_i)$, $i=1,2$, are linearly independent invariants
and it is easy to see that the subgroup of $G$ acting as the identity on $V_1$
and $V/V_1$ acts non-trivially on $\gamma_{p-1}(V_2)$ and hence the invariants
are spanned by $\gamma_{p-1}(V_1)$.

As for the second statement it is done by a computer calculation\footnote{The
source code for the calculation can be found as
\url{http://www.math.su.se/~teke/Hirokado.mg}.} though it could be done by
an unenlightening hand calculation using the restrictions obtained in the first
part.
\end{proof}
\end{proposition}
\end{subsection}
\begin{subsection}{Computation of cohomology}

The relevance of the previous section comes from the fact that by Serre duality
this duality means that the $\Gamma^*(V)$ appears as the degree $n$ cohomology
of the $\sO(m)$ on $\P^n$.
\begin{proposition}\label{DLH}
Assume that $p \ge 3$. 

\part[i] We have that $H^0(\sF,\Omega_{\sF}^1)=0$ and $H^1(\sF,\Omega_{\sF}^1)$ is $1$-dimensional
and that $H^2(\sF,\Omega_{\sF}^1)$ is $\GL_4(\F_3)$-equivariantly isomorphic to
the kernel of the map
\begin{displaymath}
\Gamma^{p-1}V\Dsum\Gamma^{2p-2}V \to
\Gamma^{p-3}V\Dsum\Gamma^{p-3}V\Dsum\Gamma^{2p-4}V\Dsum\Gamma^{p-2}V\Tensor V.
\end{displaymath}
Here $V=\Lambda^2\F_3^4$ and $\Gamma^{p-1}V$ maps to the first and last factors;
to the first factor through the contraction by the defining equation $q \in
S^2V$ for the Grassmannian and to $\Gamma^{p-2}V\Tensor V$ through the coproduct
\map{\Delta}{\Gamma^*V}{\Gamma^*V\Tensor\Gamma^*V}. Similarly, $\Gamma^{2p-2}V$
maps to the last three factors, to the first by contraction by the equation
$b(x,x^p) \in S^{p+1}V$, to the second by contraction by $q$ and to the third by
the composite of the coproduct and the $p$'th power map $\Gamma^pV \to V$ (where
we use the fact that as we are working over the prime field $\F_3$ so that $F^*V=V$).

\part[ii] We have that $H^0(\tsF,T_{\tsF})=H^0(\sF,T_{\sF})=0$.

\part[iii] We have that $H^3(\sF,\Omega_{\sF}^1)=0$.
\begin{proof}
Starting with \DHrefpart{i} everything except the explicit identification of the
map is a direct consequence of (\ref{hodge numbers}). For that identification one
goes through the proof of (\ref{hodge numbers}) and finds that
$H^2(\sF,\Omega_{\sF}^1)$ is the kernel of the map on $H^3(\P^5,-)$ induced by
the map $\sO_{\sF}(-2)\Dsum\sO_{\sF}(-p-1) \to \Omega_{\P^5|\sF}^1$ induced by
the differentials of $q$ and $b$. Furthermore, again from the proof one gets
that $H^3(\P^5,\sO_{\sF}(-2)$ is isomorphic to the kernel of
$H^5(\P^5,\sO(-p-5))\to H^5(\P^5,\sO(-4)\Dsum\sO(-p-3))$ given by multiplication
$b$ and $q$. Similarly, for $\sO_{\sF}(-p-1))$ and $\Omega_{\P^5|\sF}^1$. We
also get by the same argument that the inclusion $\Omega_{\P^5|\sF}^1 \to
\sO_{\sF}(-1)^6= \sO_{\sF}(-1)\Tensor V^*$ induces an injection on $H^3$. Now,
the composite map $\sO_{\sF}(-2)\Dsum\sO_{\sF}(-p-1) \to \Omega_{\P^5|\sF}^1 \to
\sO_{\sF}(-1)^6$ is induced on the first and second factor by the partial
derivatives of $q$ resp.\ $b(x,x^p)$. We therefore get a commutative diagram mapping
\begin{displaymath}
\begin{CD}
H^5(\P^5,\sO(-p-5))\Dsum H^5(\P^5,\sO(-2p-4))\cr
@VVV\cr
H^5(\P^5,\sO(-4))\Dsum H^5(\P^5,\sO(-p-3))\Dsum H^5(\P^5,\sO(-2p-2))\Dsum H^5(\P^5,\sO(-p-3))
\end{CD}
\end{displaymath}
to
\begin{displaymath}
\begin{CD}
H^5(\P^5,\sO(-p-3)\Tensor V^*)\cr
@VVV\cr
H^5(\P^5,\sO(-p-1)\Tensor V^*)\Dsum H^5(\P^5,\sO(-2)\Tensor V^*)
\end{CD}
\end{displaymath}
through the matrix of partial derivatives of $q$ and $b$, such that
$H^2(\sF,\Omega_{\sF}^1)$ is isomorphic to the kernel of the map induced on the
kernel of the vertical maps or, which amounts to the same the leftmost homology
of the associated double complex. Together with the identification
$H^5(\P^5,\sO(-n-6))=\Gamma^nV$ this gives the proposition except for the
identification of the maps induced by the partial derivatives. Beginning with
$q$ it is clear, using the standard covering and Cech cohomology, that for any
non-degenerate quadratic form the induced map is the composite of the coproduct
$\Gamma^nV \to \Gamma^{n-1}V\Tensor V$ and the isomorphism $V \riso V^*$ given
by the bilinear form of the form. As for $b$ it is clear that a partial
derivative of it is the $p$'th power of the same derivative of $q$.

Turning to \DHrefpart{ii} one proof is to note that as $\sF$ is normal and $\tsF
\to \sF$ has an exceptional set of codimension $2$ in $\tsF$ we have that
$H^0(\sF,T_{\sF})=H^0(\tsF,T_{\tsF})$ and then use
\cite[Prop.~2.7]{hirokado99::ae+calab+yau}. Another proof, more in the spirit of
this article is to note that, as for any complete intersection of dimension $\ge
1$, any vector field on $\sF$ is induced from one on $\P^5$. Hence the global
vector fields on $\sF$ are those global vector fields, forming the Lie algebra
$\pgl_6$, on $\P^5$ that preserve the defining ideal of $\sF$. For degree
reasons it has to preserve $q$ and is hence given by the elements of $\so_6$ or
alternatively those of $\pgl_4$ acting on $\Lambda^2V$. Now, $\pgl_4$ is
irreducible as $\GL_4(\F_p)$-representation. Indeed, by
\cite{steinberg63::repres} this is so if it is infinitesimally irreducible as
the maximal weight is the sum of two fundamental ones and for the same reason it
is infinitesimally irreducible if it is irreducible over $\GL_4$ (in
characteristic $p$) which is easily seen. Hence if there are global vector
fields, all elements of $\pgl_4$ give rise to such vector fields. This is not
true and to see this we recall the traditional form of $q$. If we use
$(q_1\co\dots\co q_6)$ as homogeneous coordinates on $\P^5$ then $q$ has the
form $q_1q_2-q_3q_4+q_5q_6$ and consequently $b(x,x^p)$ has the form
$q_1q_2^p+q_1^pq_2-q_3q_4^p-q_3^pq_4+q_5q_6^p+q_5^pq_6$. Now,
$q_1\partial/\partial q_1-q_2\partial/\partial q_2 \in \pgl_4=\so(q)$ but it
acts on $b$ to give $q_1q_2^p-q_1^pq_2$ which is easily seen not to vanish
identically on $\sF$.

Finally, \DHrefpart{iii} is a direct consequence of \DHrefpart{ii}, Serre
duality and the fact that $\omega_X$ is trivial.
\end{proof}
\end{proposition}
This result will now be combined with (\ref{resolution cohomology}). Our result will
be obtained by identifying a map in the Leray spectral sequence and then
computing its rank. The latter part is done through direct calculation and
forces us to restrict to $p=3$. We shall begin by setting up some preliminaries
that will be used to define a map that is to be identified with a differential
of the Leray spectral sequence. A further reason for the restriction to $p=3$ is
(\ref{2-invariants}), however as it seems likely that it is true for all $p$ we
shall explicitly put it in as a condition (by the same computer calculation as
in the proof of (\ref{2-invariants}) we have verified it also for $p=5,7$).
\begin{theorem}
Let \map{\pi}{\tsF}{\sF} be the natural map (over the base field $\F_p$).

\part[i] $\pi_*\Omega_{\tsF}^1=\Omega_{\sF}^1$ and
$H^0(\sF,R^1\pi_*\Omega_{\tsF}^1)$ is canonically isomorphic to 
the free $\F_p$-vector space $\F_p[\Gr(2,4)(\F_p)]$ on the set of lines of $\P^3(\F_p)$.

\part[ii] Assume that (\ref{2-invariants}) is true with $3$ replaced by
$p$. There are non-zero elements $a,b \in \F_p$ with $a+b=0$ such that the
differential
\map{d_2}{H^0(\sF,R^1\pi_*\Omega_{\tsF}^1)}{H^2(\sF,\Omega_{\sF}^1)} in the
Leray spectral sequence $E^{ij}_2=H^i(\sF,R^j\pi_*\Omega_{\tsF}^1) \implies
H^{i+j}(\tsF,\Omega_{\tsF}^1)$ is described, using \DHrefpart{i} and (\ref{DLH}), as
taking the class of an $\F_p$-rational line $\l$ to
$a\gamma_{p-1}(\l)+b\gamma_{2p-2}(\l)$.

\part[iii]\label{full hodge numbers} 
Assume $p=3$. Then the Hodge numbers of $\tsF$ are given by the
following table:
\begin{center}
\begin{tabular}{|c||c|c|c|c|}
\hline
3&1&0&0&1\\
\hline
2&0&0&42&0\\
\hline
1&0&42&0&0\\
\hline
0&1&0&0&1\\
\hline
\hline
ji&0&1&2&3\\
\hline
\end{tabular}
\end{center}
\begin{proof}
Starting with \DHrefpart{i} it follows directly from (\ref{resolution
cohomology}). As for \DHrefpart{ii} it follows from (\ref{invariants}) and
\DHrefpart{i} and the fact that $\GL_4(\F_p)$ acts on the spectral sequence that
$d_2$ must map a line to $a\gamma_{p-1}(\l)+b\gamma_{2p-2}(\l)$ for some $a$ and
$b$ that are independent of $\l$. Furthermore, by (\ref{DLH}) the element must be
mapped to zero in $\Gamma^{p-2}V\Tensor V$. Now, again by (\ref{DLH}),
$\gamma_{p-1}(f)$ is mapped to $\gamma_{p-2}(f)\tensor f$ as is
$\gamma_{2p-2}(f)$. This shows that $a+b=0$ and what remains to be shown is that
$a\ne 0$. If this were not the case then we would have that $d_2=0$. If so the
map $H^2(\sF,\Omega_{\sF}^1) \to H^2(\tsF,\Omega_{\tsF}^1)$ would be
injective. As both $\omega_{\sF}$ and $\omega_{\tsF}$ are trivial this map is
dual to $H^1(\sF,T_{\tsF}) \to H^1(\tsF,T_{\sF}^\cdot)$, where $T_{\sF}^\cdot$
is the tangent complex of $\sF$. Hence that map would be surjective.

Recall now that if \map fXY is a map of varieties then the vector space $T$ of
$\k[\epsilon]$-deformations of $f$ fits into an exact sequence
\begin{displaymath}
H^0(X,f^*T_Y^\cdot)\to T \to H^1(X,T_X^\cdot)\Dsum H^1(Y,T_Y^\cdot) \to
H^1(X,f^*T_Y^\cdot)
\end{displaymath}
(this is seen either by considering this as the problem of deforming $X$ and $Y$
and the graph of $f$ in $X\times Y$ or by using
\cite[Thm.\ III:2.1.7]{illusie71::compl} applied to the mapping cone topos of
$f$). The surjectivity of $H^1(\sF,T_{\tsF}) \to H^1(\tsF,T_{\sF}^\cdot)$
combined with the fact that $Rf_*\sO_{\tsF}=\sO_{\sF}$ would thus imply that any
first order deformation of $\sF$ would lift to a first order deformation of
$\pi$. Consider now such a deformation \map{\pi '}{\tsF'}{\sF'}. Consider a
component $E$ of the exceptional set of $\pi$. As the normal bundle, by (\ref{normal
bundle}), is isomorphic $\sO(-1)^2$ we get that $E$ lifts to a flat
$\k[\epsilon]$-curve $E'$ in $\tsF'$ and just as in the $E$-case we get that
$H^0(E',\sO_{\tsF'}/\sI_{E'}^n) \to H^0(E',\sO_{E'})=\F_p[\epsilon]$ is
surjective and as $R^1\pi_*\sO_{\tsF}=0$ we have that $\pi
'_*\sO_{\sF'}=\sO_{\sF'}$ so that by Grothendieck's formal function theorem
there is a $\k[\epsilon]$-algebra map $\sO_{\sF',x} \to \F_p[\epsilon]$. Now,
$\sO_{\sF',x}$ is a quotient of $\sO_{\P^5,x}[\epsilon]$ by an ideal $I'$ that
lifts the ideal $I$ of $\sO_{X,x}$ in $\sO_{\P^5,x}$. As $x$ is singular we have
that $I \subset \maxid^2$. The fact that we have a map $\sO_{\sF',x} \to
\F_p[\epsilon]$ then forces $I' \subset \maxid[\epsilon]$. This is not possible
for all deformations of $\sF$ as it is a complete intersection and there will
always be embedded deformations that move the ideal out of $\maxid[\epsilon]$.

Turning to \DHrefpart{iii} (and hence assuming $p=3$) we get that the kernel of
$d_2$ is equal to the kernel of the map $\F_3[\Gr(2,4)(\F_3)] \to
\Gamma^4\Lambda^2\F_3^4\Dsum \Gamma^2\Lambda^2\F_3^4$ taking $\l$ to
$\gamma_4(\l)-\gamma_2(\l)$. As $\Gamma^2V \to V\Tensor V$ given by $\gamma_2(v)
\mapsto v\tensor v$ is injective and $\gamma_4(v)$ maps to the same element,
this kernel is also the kernel of $\l \to \gamma_4(\l)$. We may describe this
map explicitly as follows:
\begin{itemize}
\item The Plücker coordinates of $\F_3$-lines are the $\F_3$-points $(q_1\co
\dots\co q_6) \in \P^5(\F_3)$ that are solutions to $q_1q_2-q_3q_4+q_5q_6=0$.

\item These are enumerated by considering separately the case $q_1=0$ and $q_1
\ne 0$.

\item In the natural basis $\gamma_{i_1}(e_1)\cdots\gamma_{i_6}(e_6)$, where
$i_1+\cdots+i_6=4$, of $\Gamma^4\F_3^6$ $\gamma_4(q_1e_1+\cdots+q_6e_6)$ has
coordinates all the degree $4$-monomials in the $q_i$.

\item This gives a $130\times 126$-matrix.

\item As a consistency check we note that for each $\F_3$-point $\l$ of $\P^3$
we get a surface $H_{\l}$ in $\tsF$ consisting of $\l$ and any line in $\P^3$
containing it. Taking its Chern class in $H^1(\tsF,\Omega_{\tsF}^1)$ and then
its image in $H^0(\sF,R^1\pi_*\Omega_{\tsF}^1)$ gives us a map $\F_3[\P^3(\F_3)]
\to \F_3[\Gr(2,4)(\F_3)]$ whose image will lie in the kernel. As the
identification $H^0(\sF,R^1\pi_*\Omega_{\tsF}^1) \riso \F_3[\Gr(2,4)(\F_3)]$ is
given by the trace map on each curve $C_P$ in $\tsF$ associated to the
$\F_3$-line $P$ we see that the map $\F_3[\P^3(\F_3)] \to \F_3[\Gr(2,4)(\F_3)]$
is given by the intersection number $(H_{\l}\cdot C_P)$. This number is easily
seen to be $1$ if $\l \subset P$ and zero otherwise.
\end{itemize}

By a computer calculation\footnote{The source code for the calculation can be
found as \url{http://www.math.su.se/~teke/Hirokado.mg}.} (which includes the
suggested consistency check) we show that the kernel of $d_2$ has dimension
$41$. As it follows from (\ref{hodge numbers}) and (\ref{DLH}) that
$h^i(\sF,\Omega_{\sF}^1)$ are equal to $0$, $1$, $89$, and $0$ respectively we
get the conclusion for $h^{1i}_{\tsF}$. For $h^{0i}_{\tsF}$ it follows from
(\ref{resolution cohomology}) and (\ref{hodge numbers}) and for the rest of the Hodge
numbers by duality.
\end{proof}
\end{theorem}
\begin{proofof}{Theorem \ref{Rigidity}}
It follows from the fact that $H^1(X,T_X)=H^1(X,\Omega_X^2)=0$ that $X$ is rigid
in characteristic $3$ so it only remains to show that it cannot be lifted to
$\W/9\W$. Assume that such a lifting $\sX$ exists. As $H^1(X,T_X)=0$ such a
lifting is unique and hence any automorphism of $X$ lifts to one of $\sX$ and as
$H^0(X,T_X)=0$ such a lifting is unique. This implies that the action of
$\GL_4(\F_3)$ on $X$ lifts to one on $\sX$.  Let now $\sL$ be the inverse image
of $\sO_{\P^5}(1)$ under the map $X = \tsF \to \sF \hookrightarrow \Gr(2,4)
\hookrightarrow \P^5$. We have that $H^1(X,\sL)=0$, indeed
$H^1(X,\sL)=H^1(\sF,\sO(1))$ by (\ref{resolution cohomology}) and $\sF$ is a complete
intersection. Now, as $H^2(X,\sO_X)=0$ $\sL$ lifts to a line bundle $\sL'$ on
$\sX$, as $H^1(X,\sO_X)=0$ such a lifting is unique. Hence, any element of
$\GL_4(\F_3)$ lifts to an automorphism of $\sL'$ compatible with the given ones
on $\sL$ and on $\sX$. The ambiguity of lifting is an element of $1+3\Z/9\Z\iso
\Z/3$ so that we get an action of a central extension of $\GL_4(\F_3)$ by
$\Z/3$. However, all such extensions are trivial by
\cite[Thm.\ 1.1]{steinberg81::gener}.

Thus we see that the natural representation of $\GL_4(\F_3)$ on $\Lambda^2\F_3$
would lift to $\W/9\W$. We shall now show that this is impossible. Indeed,
consider the element $T$ of $\GL_4(\F_3)$ that takes $e_1$ to $e_1+e_2$ and
fixes the other basis elements. It is then clear that $\Lambda^2T$ is of the
form $Id+N$ where $N \ne 0$ and $N^2=0$. We shall now show that such an element
of $\GL_6(\F_3)$ can not be lifted to an element also of order $3$ in
$\GL_6(\F_3)$. In fact, we may replace $6$ by any dimension and $3$ by any prime
$p \ge 3$. We start by taking a lifting of $N$, also denoted $N$, to $M_n(\F_p)$
with $N^2=0$. Any lifting of $Id+N$ is then of the form $Id+M$ where
$M=N+pP$. Now, we claim that $(Id+N)^p \equiv (Id+M)^p \bmod p^2$. Indeed, it is
clear that $(Id+M)^p \equiv (Id+N)^p+(N+pP)^p$. Expanding the second summand and
using that $p \ge 3$ we get in each term of the expansion always a factor $NN$
or a factor $pPpP$ both of which are zero. Now $(Id+N)^p=Id+pN$ which is never
the identity modulo $p^2$.
\end{proofof}
\end{subsection}
\end{section}
\begin{section}{de Rham-Witt cohomology}

We shall now investigate the de Rham-Witt cohomology of $X$. This should be of
independent interest but it must be admitted that our original hope was that we
should be able to show that the $n$ in (\ref{dRW}) is zero and hence that the Hodge
to de Rham spectral sequence is non-degenerate which in turn would imply that
the Hirokado threefold can not be lifted to $\Z/9$. 

For the reader's convenience we begin by a very short overview of the theory of
the cohomology of the de Rham-Witt complex (see
\cite{illusie83::finit+kuenn+rham+witt} for a more comprehensive one). The basic
objects of study are graded modules over the \Definition{Raynaud ring}, the
graded ring over the Witt ring $\W=\W(\k)$ of a perfect field $\k$ of positive
characteristic $p$, generated by $F$ and $V$ of degree zero and $d$ of degree
$1$ subject to the relations: $F\lambda = \lambda^\sigma F$, $\lambda V=V
\lambda^\sigma$ for $\lambda \in \W$ and $\sigma$ the lifting of the Frobenius
map, $FV=VF=p$, $d\lambda=\lambda d$, $d^2=0$, and $FdV=d$. A complex of
$R$-modules can then also be considered as a double complex of $\W$-modules. The
cohomology of the de Rham-Witt complex of a smooth and proper $\k$-variety $X$,
$R\Gamma(X,\W\Omega_X)$ is a complex of $R$-modules whose cohomology is
\emph{coherent}. A coherent $R$-module is  a successive extension of $R$-modules
that are either finitely generated as $\W$-module and with $d=0$ or
degree-shifts of the so called elementary dominoes, the $R$-modules of the form
\begin{displaymath}
\domU_{i}\co \pow[\k]{V} \mapright{dV^{-i}} \prod_{n \ge 0}\k dV^n
\end{displaymath}
for each integer $i$, where $FdV^0=0$ (and $VdV^n=VFdV^{n+1}=pdv^{n+1}=0$) and
we use the convention that $dV^{i}=F^{-i}d$ when $i<0$. Furthermore, there is a
right $R$-module $R_1$ such that $R_1 \Tensor_R R\Gamma(X,\W\Omega_X)$ is
canonically isomorphic to $R\Gamma(X,\Omega^\cdot_X)$, the cohomology of the de
Rham complex considered as a module over $\k[d]$ (and thus its cohomology
modules is the Hodge cohomology). In particular, as $R^i\Gamma(X,\W\Omega_X)$ is
zero if $i$ is larger than the dimension $n$ of $X$ we get
$H^n(X,\Omega^i_X)=(R_1\Tensor_RH^n(X,\W\Omega_X))$. Furthermore, the simple
complex associated to the double complex $R\Gamma(X,\W\Omega_X)$ is the complex
whose cohomology is the crystalline cohomology of $X$. The cohomology as complex
of $R$-modules of $R\Gamma(X,\W\Omega_X)$, $R^i\Gamma(X,\W\Omega_X)$ ``cuts''
the complex up into ``horisontal'' pieces which clearly differs much from the
crystalline cohomology which being the cohomology of the associated simple
complex should be thought more of as cutting up the complex along the
diagonal. There is indeed another $t$-structure (in the sense of \cite{BBD82})
on the triangulated category of bounded $R$-complexes with coherent cohomology
that cuts up this category more along the diagonal (the details are to be found
in \cite{Ek86:1}) of which we shall make extensive use.

We start with the following preliminary result.
\begin{lemma}\label{Nk}
For $k \ge 0$ there is a unique (up to isomorphism) non-trivial extension
\begin{displaymath}
\shex{\domU_{-k}[-1](1)}{N_k}{\domU_{k}}
\end{displaymath}
and we have $H^0(N)=\domU_{k-1}$ and $H^{-1}(N)=\domU_{1-k}$. In particular
$h^0(\tot(N))=h^0(\tot(N))=k-1$. Furthermore, we have
$h^{0,0}(N)=1$, $h^{1,-1}(N)=k$, $h^{1,-2}(N)=k+1$, $h^{2,-1}(N)=k+1$,
$h^{2,-2}(N)=k$, $h^{3,-3}(N)=1$, and all others are zero.
\begin{proof}
We may by shifting, $(-)(\varphi)$ (cf.\ \scite[p.~7]{Ek86:1}) assume that $k=0$
except for the statements about the Hodge numbers. Existence and uniqueness is
then \lcite[Cor.~III:1.5.4]{Ek86:1}. Consider now an $\tot$-acyclic diagonal
domino $M$. It is classified by its $F$-gauge structure $S(M)$ with
$S(M)^\infty=0$ (cf.\ \lcite[Prop.~III:1.1]{Ek86:1}) and we have that the length
of $\S(M)^i$ equals $T^{i,-i}(M)$, \lcite{Ek86:1}, (the use of $\dim_\k$ instead
of length is clearly a misprint). Consider now such an $M$ with
$T^{0,0}(M)=T^{1,-1}(M)=1$ and all other $T$'s equal to $0$. There are clearly
three corresponding $F$-gauge structures, all of them with $S(M)^0=S(M)^1=\k$
and $S(M)^i=0$ otherwise; one with $\tilde F=0,\tilde V=\symb{id}$, $\tilde
F=\symb{id},\tilde V=\symb{id}$, and $\tilde F=0,\tilde V=0$. The last case
clearly corresponds to $M\riso \domU_{0}\Dsum\domU_{0}[-1](1)$, for the second
we $\domU_{0}[-1](1)$ as a subobject with $\domU_{0}$ as quotient and the
extension is non-trivial so that $M = N_0$. In the first case we have an exact
sequence
\begin{displaymath}
\shex{\domU_{0}}{N}{\domU_{0}[-1](1)}
\end{displaymath}
and we now want to show that $H^0(N)\riso \domU_{1}$ and $H^{-1}(N)\riso
\domU_{-1}$ which implies that $N=N_1$. We have that $h^{1,0}(M)=0$ because if
it isn't then by \lcite[4.3]{Ek86:1} there would be a non-zero map $M \to
\domU_0$ which induces a non-zero map $S(M) \to S(\domU_0)$ which clearly is
impossible. Using this and considering the long exact sequence of
$H^*(R_1\Tensor^L_R -)$ applied to the above short exact sequence and using
$h^{1,i}(R_1\Tensor^L_R \domU_0)=1$, if $i=0,-2$ and $0$ otherwise
(cf.\ \cite[Cor.~3.7]{illusie83::les+rham+witt}) gives
$h^{1,0}(M)=h^{1,-1}(M)=0$. Now, by \cite[Thm.\ I:1.12]{Ek86:1}, $H^i(N)$ is
without finite torsion and is hence a domino. As $T^{0,0}(H^0(N))=T^{0,0}(N)=1$
we have that it is isomorphic to some $\domU_\ell$. Now, the fact that
$h^{1,0}(N)=h^{1,1}(N)=0$ implies that $h^{1,0}(H^0(N))=h^{1,1}(H^0(N))=0$ and
hence $\ell=1$ by \cite[Cor.~3.7]{illusie83::les+rham+witt}. On the other hand,
for the same reason $H^{-1}(N)$ is also of the form $\domU_m$ and as $N$ is
$\tot$-acyclic $H^0(\tot(H^0(N)))$ is isomorphic to $H^1(\tot(H^{-1}(N)))$ which
forces $m=-1$.

As for the Hodge numbers note to begin with that by the uniqueness $N$ is
self-dual, i.e., $D(N)=N[3](-3)$ which gives $h^{i,j}(N)=h^{3-i,3-j}(N)$. The
statement for the $h^{0,j}(N)$ is clear. For $h^{1,j}(N)$, that $H^0(N)=\domU_{k-1}$
gives $h^{1,0}(N)=0$ and then the long exact sequence of $H^i(R_1\Tensor^L_R -)$
applied to the characterising short exact sequence for $N_k$ gives the rest of
the $h^{1,j}(N)$. The rest of the Hodge number are then given by duality.
\end{proof}
\end{lemma}
We have a finite purely inseparable map from $\P^3$ blown up at the
$\F_3$-rational points to $X$. This map induces an isomorphism modulo torsion
which shows that (where we assume the base field to be $\F_3$)
$H^i(X,\W\Omega_X^j)$ is torsion except when $i=j$ and up to torsion
$H^0(X,\W)=\Z_3$, $H^1(X,\W\Omega_X^1)=\Z_3^{41}$,
$H^2(X,\W\Omega_X^2)=\Z_3^{41}$, and $H^3(X,\W\Omega_X^3)=\Z_3$ (to avoid
confusion let us recall that $\W(\F_3)=\Z_3$). Turning our attention to the
infinite part of the torsion we let as in \cite[p.~81]{illusie83::les+rham+witt}
$T^{ij}$ denote the dimension of the domino starting at
$H^j(X,\W\Omega_X^i)$. As $H^i(X,\sO_X)=0$ for $0 < i < 3$ we get that
$H^i(X,\W\sO_X)=0$ for $0 < i < 3$ and $H^3(X,\W\sO_X)$ is $V$-torsion free and as it is
torsion we have $H^3(X,\sO_X)=H^3(X,\W\sO_X)/VH^3(X,\W\sO_X)$ and its dimension
equals $T^{03}$. This gives $T^{03}=1$ and $T^{02}=0$. By duality (cf.\
\cite[3.5.1.1]{Ek84}) we have that $T^{ij}=T^{1-i,5-j}$ which gives $T^{12}=1$
and $T^{13}=0$. This gives that the only non-zero $T$'s are
$T^{03}=T^{12}=1$. Now $R^3\Gamma(X,\W\Omega_X^\cdot)$ is without finite torsion
(i.e., the torsion of the \emph{heart}, cf.\
\cite[2.4.1]{illusie83::les+rham+witt}) as it would be seen in some
$H^3(X,\Omega^j)$ (through the formula
$R_1\Tensor_RR^3\Gamma(X,\W\Omega_X^\cdot)=H^3(X,\Omega^*_X)$ of
\cite[Thm.\ II:1.2]{illusie83::les+rham+witt}) and all of that cohomology is
already accounted for by the fact that $T^{03}=1$ and
$H^3(X,\W\Omega_X^3)=\Z_3$. By duality, \cite[Thm.\ IV:3.5]{Ek84}, there is no
nilpotent torsion in any $R^*\Gamma(X,\W\Omega_X^\cdot)$ and as also
$H^2(X,\W\sO_X)=0$ we see that the only possibilities for semi-simple torsion is
in $H^2(X,\W\Omega_X^j)$ for $j=1,2$. This completes the preliminary results on
the way to proving the following theorem. In it we shall for simplicity work over
an algebraic closure of $\F_3$ even though the result can be shown (e.g.,
\emph{a posteriori} by Galois descent) to be true over $\F_3$.
\begin{theorem}\label{dRW}
\part The de Rham-Witt cohomology over $\overline{\F}_3$ of the Hirokado threefold $X$
is isomorphic to the direct sum (as an object in the derived category of graded
modules over the Raynaud ring $R$)
\begin{displaymath}
\W\Dsum\W^{41}[-1](-1)\Dsum\W^{41}[-2](-2)\Dsum\W[-3](-3)\Dsum M[-3]
\end{displaymath}
where $\W:=\W(\overline{\F}_3)$.
Furthermore, there is an $n \ge 0$ such that $M$ fits into a distinguished
triangle
\begin{displaymath}
 \to \W/p^n\W[1](-1) \to M \to N' \to  \W/p^n\W[2](-1) \to
\end{displaymath}
and $N'$ fits into a distinguished triangle
\begin{displaymath}
\to N \to N' \to \W/p^n\W[-1](-1) \to N[1] \to 
\end{displaymath}
where $N$ is a diagonal domino (cf.\
\cite[Def.~III:2.1]{Ek86:1}) which is the unique (up to isomorphism) extension
of $\domU_{-1}[1](-1)$ by $\domU_1$ (cf.\
\cite[I:2.4.13]{illusie83::les+rham+witt}). It is self-dual (up to an
appropriate shift), the map $ N' \to  \W/p^n\W(-1)$ is determined by the
composite with $N \to N'$ and that composite is dual to $\W/p^n\W[-1](-1) \to
N[-1]$. Finally, the map $N \to \W/p^n\W(-1)$ is the unique (up to isomorphism)
non-zero such map.

\part The crystalline cohomology of $X$ has the following values:
$H^0(X/\W)=H^6(X/\W)=\W$,
$H^1(X/\W)=H^5(X/\W)=0$, $H^2(X/\W)=\W^{41}$, $H^3(X/\W)=\W/p^n\W$, and
$H^4(X/\W)=\W/p^n\W\Dsum\W^{41}$ (with the same $n$ as in the previous
part). The Hodge to de Rham spectral sequence is degenerate precisely when
$n>0$.
\begin{proof}
We have that $R^0\Gamma(X,\W\Omega_X^\cdot)$ is always torsion free so the
computation modulo torsion gives $R^0\Gamma(X,\W\Omega_X^\cdot)=\W$. Similarly,
using that $T^{03}$ and $T^{12}$ are the only non-zero $T$'s and that
$R^1\Gamma(X,\W\Omega_X^\cdot)$ contains no finite torsion we get that it also
is torsion free and so $R^1\Gamma(X,\W\Omega_X^\cdot)=\W^{41}(-1)$. By duality
the map $\W = R^0\Gamma(X,\W\Omega_X^\cdot) \to R\Gamma(X,\W\Omega_X^\cdot)$
gives a map $R\Gamma(X,\W\Omega_X^\cdot)\to \W[-3](-3)$ (which is just the trace
map) and from \cite[Cor.~III:1.5.4~i]{Ek85} this map splits (it can also be
split directly by considering the cohomology class of a point). This allows us
to split off as a factor $\W[-3](-3)$ as well . Let
$R^{[1,5]}\Gamma(X,\W\Omega_X^\cdot)$ be the remaining complex. This time we
have a map $\W^{41}[-1](-1) \to R^{[1,5]}\Gamma(X,\W\Omega_X^\cdot)$ and a dual
map $R^{[1,5]}\Gamma(X,\W\Omega_X^\cdot) \to \W^{41}[-2](-2)$. Again using
\cite[Cor.~III:1.5.4 i]{Ek85} this map has a section and we can split off
$\W^{41}[-2](-2)\Dsum\W^{41}[-3](3)$. We are left with a torsion complex of the
form $M[-3]$ and we can use the previously obtained results together with (\ref{full
hodge numbers}) to get that
\begin{itemize}
\item The Hodge numbers of $M$ are
$h^{00}(M)=h^{2,-1}(M)=h^{1,-2}(M)=h^{3,-3}(M)=1$ and the rest are zero.

\item $T^{00}(M)=T^{1,-1}=1$ and the rest are zero.

\item Only $H^{-1}(M)^{1}$ and $H^{-1}(M)^{2}$ can contain finite torsion and it
is semi-simple.

\item We have a duality isomorphism $D(M)=M[-3](3)$, where $D(-)$ is defined in
\cite[Def.~III:2.8]{Ek84}.
\end{itemize}
We shall now use these properties to describe $M$. We shall consider the
diagonal cohomology (cf.\ \scite[Def.~I:1.2]{Ek86:1}) of $M$. By the properties
of $M$ we get that $\tH^i(M)=0$ if $i \ne 0,1$ and that $\tH^1(M)$ is the finite
torsion part of $H^{-1}(M)^{2}$. On the other hand we have $t^1(\tH^0(M))$, the
maximal finite torsion subobject of $\tH^0(M)$ (cf.\ \lcite[I:1.9]{Ek86:1}) and
the duality hypothesis on $M$ implies (cf.\ \lcite[I:Prop.~1.11]{Ek86:1}) that
$\tH^1(M)$ is dual to $t^1(\tH^0(M))$ and that $N:=\tH^0(M)/t^1(\tH^0(M))$ is
self-dual. Consider now the distinguished triangle
\begin{displaymath}
\to t^1(\tH^0(M)) \to M \to M' \to t^1(\tH^0(M))[1] \to
\end{displaymath}
and the associated long exact sequence of Hodge cohomology in degree $1$.
Using the Hodge numbers of $M$ gives the long exact sequence
\begin{displaymath}
0 \to {}_3T \to \oF_3 \to H^{-2}(R_1\Tensor_RM')^1 \to T/3T \to 0
\end{displaymath}
where $T := H^{-1}(t^1(\tH^0(M)))^1$ and that $ H^{i}(R_1\Tensor_RM')^1$ for $i
\ne -2$. The exact sequence gives that the dimension of ${}_3T$ is $1$ or $0$
and as $T$ is torsion and finitely generated as $\W$-module that means that
$T=\W/p^n\W$ for some $n \ge 0$. The exact sequence then also gives that
independently of the value of $n$ we have $h^{1,-2}(M')=1$ and as $\tH^1(M)$ is
concentrated in degree $2$ we have $h^{1,i}(M')=h^{1,i}(N)$. Note also that as $
H^{-2}(R_1\Tensor_RM')^1 \to T/3T$ is non-zero so is $M' \to t^1(\tH^0(M))[1]$.
Now, as $N$ is a diagonal torsion complex without finite torsion with
$T^{0,0}=T^{1,-1}=1$ we get that $H^0(N)$ and $H^{-1}(N)(1)$ are $1$-dimensional
dominoes as they are without finite torsion (cf.\ \lcite[Thm.\ 1.12]{Ek86:1}) and
hence isomorphic to $\domU_k$ and $\domU_m$ for some $k$ and $m$ (cf.\
\cite[Prop.~I:2.15]{illusie83::les+rham+witt}) and by auto-duality for $N$ we
have $m=-k$. Furthermore, $0=\tH^0(R_1\Tensor_RN)^1=\tH^0(R_1\Tensor_RH^0(N))^1$
and so by \cite[I:2.14.6]{illusie83::les+rham+witt} $k > 0$. Considering now the
distinguished triangle
\begin{displaymath}
\to \domU_{-k}[1](-1) \to N \to \domU_k \to \domU_{-k}[2](-1)
\end{displaymath}
and the associated long exact sequence together with
\cite[Cor.~I:3.7]{illusie83::les+rham+witt} gives
\begin{displaymath}
0  \to \oF_3 \to \oF_3^{k+1} \to \oF_3 \to 0 \to \oF_3^{k-1} \to 0
\end{displaymath}
which gives $k=1$. Note here also that as $\oF_3^{k+1} \to \oF_3$ is non-zero so
is the map $\domU_k \to \domU_{-k}[2](-1)$. Hence by Lemma \ref{Nk} $N$ is
isomorphic to $N_k=N_1$.

Similarly, one shows that up to isomorphism there is only one non-zero map $N
\to t^1(\tH^0(M))[1]$ and that $\tH^0(M)$ is unique up to isomorphism. Then the
map $\tH^1(M) \to N[2]$ is dual to this map and one shows easily that $\tH^0(M)
\to N$ induces a bijection between maps $\tH^1(M) \to N[2]$ and maps $\tH^1(M)
\to \tH^0(M)[2]$.

Turning to the second part of the theorem we note that by Lemma \ref{Nk} $N\riso
Nk$ is $\cs$-acyclic (cf.\ \cite[Def.~I:3.1]{Ek86:1}). The exact sequence
\begin{displaymath}
0 \to T \to \tH^0(M) \to N \to 0
\end{displaymath}
now gives $H^0(\cs(\tH^0(M)))=T$ and $H^1(\cs(\tH^0(M)))=0$ and then
$H^0(\cs(M))=T$ and $H^1(\cs(M))=H^0(\cs(\tH^1(M)))=T^*$, the latter by
duality. This gives the de Rham cohomology of $X$ and it is also clear that the
Hodge to de Rham spectral sequence degenerates precisely when crystalline
cohomology has torsion.
\end{proof}
\end{theorem}
\begin{remark}
It follows from \cite{deligne87::relev+rham} that if $n=0$, then $X$ does not lift to
$\Z/9$. Unfortunately we do not know the value of $n$ and have therefore been
forced to prove non-liftability in a different fashion.
\end{remark}
\end{section}
\begin{section}{Smooth K3-pencils on Calabi-Yau threefolds}
\label{sec:pencils}

In \scite{schroeer03::some+calab+yau+witt} examples of non-liftable Calabi-Yau
threefolds were constructed as smooth pencils of K3-surfaces. In this section we
shall study such pencils. Note that by a \Definition{supersingular K3-surface}
we shall mean, unless otherwise mentioned, a K3-surface for which the rank of
the N\'eron-Severi group is equal to $22$.

We start by recalling a construction from
\lcite{schroeer03::some+calab+yau+witt} as well as correcting what we claim is a
mistake in it. In $\Char\k=2$, \cite{schroeer03::some+calab+yau+witt} constructs
a family of generalised Kummer surfaces associated to the Moret-Bailly family and
an automorphism of order $3$ of each fibre, the family being the minimal
resolution of the quotient which fibre by fibre is the minimal resolution of the
quotient. The Artin invariant of these resolutions have been computed in
\cite{schroeer03::some+calab+yau+witt} but we give here an alternative way of
performing this calculation.
\begin{proposition}
Let $E$ be a supersingular elliptic curve in characteristic $2$ and $\zeta$ an
automorphism of it of order $3$ fixing the origin of $E$. Let $\alpha_2 \subset
E\times E$ be an embedding of group schemes and $A := E\times E/\alpha_2$ and
denote again by $\zeta$ the automorphism of $A$ induced by the diagonal action
of $\zeta$ on $E\times E$. Then the minimal resolution $X$ of $A/\zeta$ is a
supersingular K3-surface with Artin invariant $\sigma_0$ equal to $1$ if $A$ is
isomorphic to $E\times E$ and equal to $2$ if not.
\begin{proof}
That it is supersingular follows directly for instance from the fact that the
rank of the N\'eron-Severi group of $A$ equals $b_2$ and which means that the
Brauer group of $A$ is finite which implies that of $X$ also is. Hence we
have to show that the discriminant of the N\'eron-Severi group, $NS$, of $X$ is
$4$ or $16$. For that we may localise it at $2$ which we shall indeed do. Having
done that we have that $NS(A/\zeta)$ equals the fixed point subgroup of $NS(A)$
and as the cup product is multiplied by the degree, which is $3$, we see that
they have the same discriminant (at $2$). As all the singularities of $A/\zeta$
are rational double of type $A_2$ (cf.\ \cite[Prop.\
5.2]{schroeer03::some+calab+yau+witt}) we have (outside of $3$) $NS =
NS(A/\zeta)\perp M$, where $M$ is spanned by the exceptional curves for the
resolution map $X \to A/\zeta$. Now, again as the singularities are of type
$A_2$ the discriminant of $M$ is a power of $3$ so that the discriminant of $NS$
(still at $2$) equals that of $NS(A/\zeta)$ and hence that of $NS(A)^\zeta$. We
shall now follow \scite[Prop.\ 6.9]{Og79} in determining the discriminant of
$NS(A)^\zeta$. Now (cf.\ \lcite[Cor.\ 1.6]{Og79}), $NS(A)\Tensor\Z_2$ equals
under the first Chern class map the submodule of $H^2(A/\W)$ of elements $w$ for
which $Fw=2w$. Furthermore, as $3$ and $2$ are relatively prime we have that
$H^2(A/\W)$ splits up under the action of $\zeta$ as a direct (orthogonal under
the cup product) sum of $F$-crystals of the invariant part $H^2(A/\W)^\zeta$ and
its complement $M$. We know that $H^2(A/\W)$ divided by the sub-$\W$-module of
elements $w$ for which $Fw \in 2H^2(A/\W)$ is canonically isomorphic to
$H^2(A,\cO_A)$. The action of $\zeta$ on $H^2(A,\cO_A)$ is trivial; one may for
instance use duality and the fact that it acts trivially on non-trivial
$2$-forms as the determinant of the action on the tangent space of a fixed point
is $1$ (cf.\ \cite[proof of Prop.\ 5.2]{schroeer03::some+calab+yau+witt}). Hence
$F$ is divisible by $2$ on $M$ and the cup product on $M^{F=2}$ is a perfect
pairing. Thence the discriminant of $H^2(A/\W)^\zeta$ equals that of $H^2(A/\W)$
so it is enough to compute the latter.  For that we start by making computations
for $E\times E$. We have that $H^1(E/W)$ is isomorphic to the $F$-crystal
generated by an element $u$ with the relation $F^2u=-2u$. In particular $u$ and
$Fu$ form a $\W$-basis for $H^1(E/W)$. Then $H^1(E\times E/\W)$ is the direct
sum of two copies of $H^1(E/\W)$ whose generators will be denoted $u$ and $v$
and thus a $\W$-basis is given by $u,Fu,v,Fv$. This means that $H^2(E\times
E/\W)$ has a basis consisting of $u\wedge v,Fu\wedge Fv,u\wedge Fv,Fu\wedge
v,u\wedge Fu,v\wedge Fv$. It is then easy to see that the elements of the kernel
of $F-2$ can uniquely be written in the form $au\wedge Fu+bv\wedge Fv+2cu\wedge
v+c^\sigma Fu\wedge Fv+du\wedge Fv-d^\sigma Fu\wedge v$, where $a,b \in \Z_2$,
$c,d \in \W(\F_4)$, and $\sigma$ is the Frobenius map on $\W(\F_4)$. In
particular we see that they generate a sub-$\W$-module of colength one in $H^2(E
\times E/\W)$ so that the discriminant of the kernel of $F-2$ is $4$. Let now
$(s\co t) \in \P^1(\k)$ and consider the sub-$\W$-module of $H^1(E\times E/\W)$
spanned by $Fu$, $Fv$, $pu$, $pv$, and $w := \underline{s}u+\underline{t}v$,
where $\underline{r} \in \W(\k)$ is the Teichm\"uller representative given in
vector form by $(r,0,0,\dots)$. This is a sub-$F$-crystal and is the Dieudonn\'e
module for a typical member of the Moret-Bailly pencil so that we may assume
that corresponds to $A$. We may also assume that $st \ne 0$ because if not then
the member is the product of supersingular elliptic curves and then we may also
put $t=1$. We then have that $w\wedge Fu,w\wedge Fv,w\wedge 2u,Fu\wedge
Fv,Fu\wedge 2u,Fv\wedge 2u$ is a $\W$-basis for $H^2(A/\W)$. Writing them out
gives $\underline{s}u\wedge Fu+v\wedge Fu,\underline{s}u\wedge Fv+v\wedge Fv,2
u\wedge v,Fu\wedge Fv,2Fu\wedge u,2Fv\wedge u$ and we may change the first to
$u\wedge Fu+\underline{s^{-1}}v\wedge Fu$. That means that the kernel of $F-2$
contains $apu\wedge Fu+bpv\wedge Fv+pcu\wedge v+c^\sigma Fu\wedge Fv+du\wedge
Fv-d^\sigma Fu\wedge v$, with $a,b,c,d$ as above and it is larger precisely if
there are $e,f \in \W$ not both divisible by $2$ such that $a(u\wedge
Fu+\underline{s^{-1}}v\wedge Fu)+b(\underline{s}u\wedge Fv+v\wedge Fv)$ lies in
the kernel of $F-2$. This first forces $a,b \in \Z_2$ and we get
$a\underline{s^{-1}}=(b\underline{s})^\sigma$ which gives $a=b\underline{s^3}$
and as $s \ne 0$ and not both of $a$ and $b$ are divisible by $p$ this implies
$s^3 \in \F_2$, i.e., $s \in \F_4$. Conversely, if $s \in \F_4$ we can reverse
the steps so that the discriminant equals $4$ if $s \in \F_4$ and $16$ if
not. Now, the automorphism group of $H^1(E\times E/\W)$ acts transitively on
$\P^1(\F_4)$ so we see that $s \in \F_4$ precisely when $A \riso E\times E$.
\end{proof}
\end{proposition}
We recall the theory of the period map for supersingular K3-surfaces (cf.\
\cite{Og79} and \cite{Og82}). Recall to start with that if $T$ is non-degenerate
even lattice, then its \Definition{discriminant group}, $D(T)$, is the quotient
$T^*/T$ where the embedding $T \to T^*$ is given by the scalar product. On
$D(T)$ we may define a non-degenerate $\Q/\Z$-valued quadratic form by $\psi(x)
= x^2/2 \bmod \Z$. A \Definition{K3-lattice of invariant $\sigma_0$}, $1 \le
\sigma_0 \le 10$ is a lattice $T$ of index $(1,21)$ and for which $D(T)\riso
(\Z/p)^{2\sigma_0}$ with a non-split $1/p\Z/\Z$-valued quadratic form. Such a
lattice is unique up to isometry and we shall fix a choice of one such
lattice. We let $T_0 \subseteq T\Tensor\Z/p$ be the radical of the induced
quadratic form. Multiplication by $p$ induces an isomorphism $D(T) \riso T_0$
and we use it to transfer the quadratic form to $T_0$. We also let $T' \subset
T\Tensor\Z/p$ be a complement to $T_0$ so that in particular the restriction of
$\psi$ to $T'$ is non-degenerate.

Let now $\k$ be a field of positive characteristic $p$ and \map{\pi}{X}{S} be an
arbitrary family of supersingular K3-surfaces with $S$ smooth over $\k$. We
assume given a \Definition{$T$-marking} of $\pi$, i.e., an isometric embedding
$T\times S \to \underline{\Pic}(X/S)$.  The first Chern class gives a map
$T\Tensor\cO_{S} \to R^1\pi_*\Omega^\cdot_{X/S}$ and, as proved in \scite{Og79},
the kernel, which necessarily is contained in $T_0\Tensor\cO_{S}$, is of the
form $F^*_{S}K$ for some sub-bundle $K \subset T_0\Tensor\cO_{S}$ that is
maximal totally isotropic with respect to the discriminant pairing on $T_0$,
hence of rank $\sigma_0$ and for which $K \cap F^*K$ is of rank
$\sigma_0-1$. Furthermore, \lcite{Og79}, the image of the first Chern class map
is contained in $F^2R^1\pi_*\Omega^\cdot_{X/S}$ of the Hodge filtration and the
kernel of the composite $T\Tensor\cO_{S} \to F^2R^1\pi_*\Omega^\cdot_{X/S} \to
F^2/F^1=R^1\pi_*\Omega^1_{X/S}$ is $K+F^*K$. If we let $\cM^T$ be the moduli
space of such subspaces of $T_0$ then we thus get a map $S \to
\cM^T_{\sigma_0}$. By results of Ogus and Rudakov-Shafarevich (cf.\ \cite{Og82}
and \cite{rudakov81::surfac+k3}) is essentially a moduli space for marked
K3-surfaces once we admit Proposition \ref{autom free}. There are some
subtleties but a consequence is that there exists a family of $T$-marked
K3-surfaces over $\cM^T$ such that its pullback along the map $S \to \cM^T$ is
isomorphic to the given family $X \to S$ (though not necessarily as
\emph{marked} surfaces).
\begin{proposition}\label{autom free}
An automorphism acting trivially on the Picard group of a supersingular
K3-surface is the identity.
\begin{proof}
According to \cite[Cor.\ 2.5]{Og79} the only possible counter example is a
superspecial K3-surface in characteristic $2$. However, according to
\cite{rudakov78::super+k3} there is a unique such supersingular (in the sense of
$\rho=22$) surface and it is then the generalised Kummer surface associated to
the square of the supersingular elliptic curve and an automorphism of order $3$
(cf.\ \cite[Prop. 5.1]{schroeer03::some+calab+yau+witt}). As both the curve and
the automorphism lifts over $\W(\k)$ so does the K3-surface and then the proof
of \cite[Cor.\ 2.5]{Og79} goes through.
\end{proof}
\end{proposition}
We have a diagram, \lcite[5.4.1]{Og79},
\par\noindent\parbox{0.9\textwidth}{\begin{displaymath}
\begin{xy}
\xymatrix{
0\ar[r]&F^*K\ar[r]&T\Tensor\cO_{S}\ar[r]&R^2\pi_*\Omega^\cdot_{X/S}\ar[r]&F^*K\ar[r]&0\\
0\ar[r]&F^*K\ar[r]\ar@2{-}[u]&T\Tensor\cO_{S}\ar@2{-}[u]\ar[r]&
   F^1R^2\pi_*\Omega^\cdot_{X/S}\ar[r]\ar@{^{(}->}[u]&K \cap F^*K\ar@{^{(}->}[u]\ar[r]&0\\
0\ar[r]&F^*K\ar[r]\ar@2{-}[u]&K+F^*K\ar@{^{(}->}[u]\ar[r]&
   F^2R^2\pi_*\Omega^\cdot_{X/S}\ar[r]\ar@{^{(}->}[u]&0\ar@{^{(}->}[u].
}
\end{xy}
\end{displaymath}
}\hfill\parbox{0.05\textwidth}{\begin{equation}\label{fundamental diagram}%
\end{equation}}\par
Put $L := K/K \cap F^*K$, $L' := F^*K/K \cap F^*K$, and $M := \det(K \cap
F^*K)^{-1}$. We have that as $K \cap F^*K$ and $K+F^*K$ are subbundles that
$K+F^*K/K \cap F^*K$ is the direct sum of $L$ and $L'$ and as $K+F^*K$ are
annihilators of each other we have that $L\Tensor L' = \det(K+F^*K/K \cap
F^*K)=\det(T_0\Tensor\cO_{S})=\cO_{S}$ and thus $L' = L^{\tensor-1}$. On
the other hand $\det(K)=L\Tensor M^{\tensor-1}$ and thus $\det(F^*K)=L^p\Tensor
M^{\tensor-p}$ but on the other hand $\det(F^*K)=L'\Tensor
M^{\tensor-1}=L^{\tensor-1}\Tensor M^{\tensor-1}$ which together gives
$L^{\tensor p+1}=M^{\tensor p-1}$. On the other hand, 
It follows from (\ref{fundamental diagram}) that
\begin{displaymath}
\cO_{S}(-2)\riso
R^2\pi_*\cO_{X/S}=R^2\pi_*\Omega^\cdot_{X/S}/F^1R^2\pi_*\Omega^\cdot_{X/S}\riso F^*K/K\cap F^*K\riso
L^{-1}
\end{displaymath}
and $\det K = L\Tensor M^{\tensor-1}$ so that we get
\begin{equation}\label{line equation}
(\det K)^{\tensor p-1} \riso (R^2\pi_*\cO_X)^{\tensor2}
\end{equation}
Note that in particular on $\cM^T$, as $K$ is the tautological
subbundle $\det K$ is anti-ample, and so is $R^2\pi_*\cO_X$ or by duality
$\pi_*\omega_{X/\cN^T_{\sigma_0}}$.

Let now \map{\pi}{X}{\P^1} be a smooth pencil of supersingular (in the sense of
$\rho=22$) K3-surfaces on a Calabi-Yau threefold over a field $\k$ of
characteristic $p$. Let $\sigma_0 = \sigma_0(\pi)$ be the Artin invariant of its
generic fibre. Let also $T$ be a fixed supersingular K3-lattice of discriminant
$-p^{2\sigma_0}$, let $T_0$ be the radical of $T\Tensor\Z/p$ and $T' \subseteq
T\Tensor\Z/p$ a complement of $T_0$. We start by proving some results on the
cohomology of $X$ and $\pi$, some of which are already contained in
\scite{schroeer03::some+calab+yau+witt} but as \lcite[Lemma
1.4]{schroeer03::some+calab+yau+witt} has an incorrect proof (notably the first
statement of that proof is definitely false) we shall provide proofs for them where
appropriate.
\begin{proposition}\label{prel cohomology}

\part[i] $R^i\pi_*\Z_\ell$ are constant \'etale sheaves of rank $1$, $0$, $22$,
$0$, and $1$ for $i=0,1,2,3,4$ respectively, where $\ell$ is a prime different
from $p$.

\part[ii] We have $b_1(X)=b_3(X)=b_5(X)=0$, $b_2(X)=b_4(X)=23$, and
$b_0(X)=b_4(X)=1$.

\part[iiia] We have $R^2\pi_*\cO_X\riso \cO_{\P^1}(-2)$.

\part[iii] $\underline{\Pic}(X/\P^1)\Tensor\Z[1/p]$ is a constant sheaf (of rank
$22$). In particular, if $\overline{\eta}$ is the generic point of $\P^1$,
$\Gamma(\P^1,\underline{\Pic}(X))=\Pic(X_{\overline{\eta}})$ for which we may (and
shall) choose an isomorphism with $T$.

\part[iv] $h^{0i}(X)=h^{3i}(X)=0$ for $i=1,2$ and $h^{0i}(X)=h^{3i}(X)=1$ for $i=0,3$.

\part[v] $-\chi(\Omega^1_X)=\chi(\Omega^2_X)=24$.

\part[vi]\label{DR cohomology} The relative de Rham-cohomology sheaf
$H^2_{DR}(X/\P^1)=R^2\pi_*\Omega^\cdot_{X/\P^1}$ is locally free and hence has
an integrable connection $\nabla$, the Gauss-Manin connection. The Chern class
map $T \to H^2_{DR}(X/\P^1)$ maps $T'\Tensor\cO_{\P^1}$ to a non-degenerate
horisontal submodule. Let $H'$ be its orthogonal complement. The dimension of
$H^2_{DR}(X/\k)$ equals $23-2\sigma_0$ plus the dimension of the space of global
horisontal sections of $H'$.
\begin{proof}
The first part follows from the proper and smooth base change theorems and the
second part is an immediate consequence of the first. For \DHrefpart{iiia} we
have $\cO_X\riso\omega_X=\pi^*\omega_{\P^1}\Tensor\omega_{X/\P^1}$ and thus
$R^2\pi_*\co_X\riso \omega_{\P^1}\Tensor R^2\pi_*\omega_{X/\P^1}$ and by
relative duality $R^2\pi_*\omega_{X/\P^1} = \cO_{\P^1}$.

As for \DHrefpart{iii} it follows from \DHrefpart{i} as
$\underline{\Pic}(X/\P^1)\Tensor\Z_\ell \to R^2\pi_*\Z_\ell$ is an
isomorphism. Continuing with \DHrefpart{iv}, it is \lcite[Prop.\
1.1]{schroeer03::some+calab+yau+witt} and duality. Now, \DHrefpart{v} follows
from the fact that $\chi(\Omega^1_X)=-\chi(\Omega^2_X)$ by duality, that
$\chi(\cO_X)=\chi(\Omega^3_X)=0$ and the fact the alternating sum of the
$\chi$'s equal to the Euler characteristic which is $48$ by \DHrefpart{ii}.

Turning to \DHrefpart{vi}, the local freeness and the rank statement is clear as
the de Rham cohomology of K3-surface is independent of the surface. As the
cup-product on $T'$ is an non-degenerate we get a splitting
$H^2_{DR}(X/\P^1)=T'\Tensor\cO_{\P^1}\Dsum H'$ which is preserved by the
connection as the cup product and the Chern classes are horisontal. Furthermore,
the connection on $T'\Tensor\cO_{\P^1}$ is the flat one and $T'$ is of dimension
$22-2\sigma_0$. Hence the statement is equivalent to saying that the dimension
of $H^2_{DR}(X/\k)$ equals $1$ plus the dimension of the space of flat global sections
of $H^2_{DR}(X/\P^1)$. The wedge filtration gives a short exact sequence of complexes.
\begin{displaymath}
\begin{CD}
\cO_X @>d>> \Omega^1_{X/\P^1}  @>d>> \Omega^2_{X/\P^1}\\
      @|                  @AAA    @AAA           \\
\cO_X @>d>> \Omega^1_{X}       @>d>> \Omega^2_{X}  @>d>> \Omega^3_{X}\\
      @.                   @AAA           @AAA               @AAA    \\
      @.    \pi^*\Omega^1_{\P^1} @>d>> \pi^*\Omega^1_{\P^1}\Tensor\Omega^1_{X/\P^1} @>d>> \pi^*\Omega^1_{\P^1}\Tensor\Omega^2_{X/\P^1}
\end{CD}
\end{displaymath}
As we know that $H^i_{DR}(X/\P^1)=0$ for $i\ne 0,2,4$ we get exact sequences
\begin{displaymath}
\begin{array}{@{0\to}ccccccc@{\to0}}
   \pi_*\Omega^\cdot_X&\mapright{}&H^0_{DR}(X/\P^1)&\mapright{\nabla}&
   \Omega^1_{\P^1}\Tensor H^0_{DR}(X/\P^1)&\mapright{}&R^1\pi_*\Omega^\cdot_X\\
R^2\pi_*\Omega^\cdot_X&\mapright{}&H^2_{DR}(X/\P^1)&\mapright{\nabla}&
   \Omega^1_{\P^1}\Tensor H^2_{DR}(X/\P^1)&\mapright{}&R^3\pi_*\Omega^\cdot_X\\
R^4\pi_*\Omega^\cdot_X&\mapright{}&H^4_{DR}(X/\P^1)&\mapright{\nabla}&
   \Omega^1_{\P^1}\Tensor H^4_{DR}(X/\P^1)&\mapright{}&R^5\pi_*\Omega^\cdot_X,\\
\end{array}
\end{displaymath}
where we have used \cite[3.2.5]{Ka70} to identify the boundary maps with the
connections on the de Rham cohomology. We have canonical identifications
$H^0_{DR}(X/\P^1)=H^2_{DR}(X/\P^1)=\cO_{\P^1}$ with the flat connection and in
particular, using the Cartier isomorphism and the first exact sequence, we get
that $R^1\pi_*\Omega^\cdot_X\riso\Omega^1_{\P^1}$. The spectral sequence
$H^i(\P^1,R^j\pi_*\Omega^\cdot_X) \implies H^{i+j}(X,\Omega^\cdot_X)$
degenerates for trivial reasons and hence gives a short exact sequence.
\begin{displaymath}
\shex{H^1(\P^1,R^1\pi_*\Omega^\cdot_X)}{H^2_{DR}(X/\k)}{H^0(\P^1,R^2\pi_*\Omega^\cdot_X)}
\end{displaymath}
and by the second exact sequence above we have that
$H^0(\P^1,R^2\pi_*\Omega^\cdot_X)$ is the space of global horisontal sections of
$H^2_{DR}(X/\P^1)$. On the other hand,
$H^1(\P^1,R^1\pi_*\Omega^\cdot_X)=H^1(\P^1,\Omega^1_{\P^1})=\k$. This gives
\DHrefpart{vi}.
\end{proof}
\end{proposition}
\begin{remark}
If we drop the assumption that the fibres are supersingular then we still get
that $R^2\pi_*\Z_\ell(1)$ is constant and also that
$H^2(X,\Z_\ell(1))=H^2(\P^1,\Z_\ell(1))\Dsum H^0(\P^1,\Z_\ell(1))$. On the other
hand, as $h^{02}=0$ we get that all the eigenvalues of Frobenius (assuming that
we specialise for the moment to a finite field $\F_q$) are divisible by $q$ and
hence are of the form $q$ times a root of unity. By the triviality of
$R^2\pi_*\Z_\ell(1)$ the same is then true for all the (closed) fibres of
$\pi$. This shows that the fibres of $\pi$ are supersingular in the sense of
the height of the Brauer group being infinite. Under the Tate conjecture it is
then also supersingular in the sense of $\rho=22$.
\end{remark}
From now on we shall choose an isometry with $\Pic(X_{\overline{\eta}})$ and $T$
and thus $\pi$ becomes a family of $T$-marked surfaces. We shall also follow the
proposition and choose a complement $T'$ to $T_0 \subset T\Tensor\Z/p$ and let
$H'$ be the orthogonal complement of the image of $T'\Tensor\cO_{\P^1}$ in
$H^2_{DR}(X/\P^1)$. Similarly, we let $H''$ be the orthogonal complement of the
image of $T'\Tensor\cO_{\P^1}$ in $R^1\pi_*\Omega^1_{X/\P^1}$.
\begin{proposition}\label{quotient possibilities}

\part[KS] $p \le 3$ and if the Kodaira-Spencer map of $\pi$ is zero then $p=2$.

\part[det] $K/K\cap F^*K \riso \cO_{\P^1}(2)$ and $\det K\riso\cO_{\P^1}(m)$,
where $m=-6$ if $p=2$ and $m=-4$ if $p=3$. 

\part[KS1] If the Kodaira-Spencer map of $\pi$ is non-zero then there is a
non-zero map $K \cap F^*K \to \cO_{\P^1}(-4)$.

\part[van] $H^0(\P^1,K)=0$ and $\sigma_0 \le 4$ if
$p=2$ and $\sigma_0=2$ if $p=3$.

\part[ii] If the Kodaira-Spencer map is non-zero, then $h^{1i}(X)=h^{1i}(X)=0$
for $i=0,3$ and $h^{11}(X)$ and $h^{12}(X)$ have values according to the
following table:
\begin{displaymath}
\begin{array}{|c|c||c|c|}
\hline
p&\sigma_0&h^{11}&h^{12}\\
\hline
2&2&26&2\\
2&3&25&1\\
2&4&24&0\\
3&2&24&0\\
\hline
\end{array}
\end{displaymath}

\part[iii]\label{DR cohomology1} $\dim_{\k}H^2_{DR}(X/\k)=\delta+23-\sigma_0$,
where $\delta=4$ if $p=2$ and $\delta=2$ if $p=3$.
\begin{proof}
If the Kodaira-Spencer map is zero then the map from $\P^1$ to the space of
marked K3-surfaces factors through the Frobenius map on $\P^1$ and hence
$R^1\pi_*\cO_X$ is the pullback by of some line bundle by the Frobenius map and
in particular its degree is divisible by $p$. Proposition \ref{prel cohomology}
then shows that $p=2$.

Now, putting $\det K=\cO_{\P^1}(-m)$ and using (\ref{line equation}) and
Proposition \ref{prel cohomology} we get $4=(p-1)m$ and this
implies that $p \le 5$ and $(p,m)$ is one of the pairs $(2,4)$, $(3,2)$ or
$(5,1)$. This proves \DHrefpart{det}.

Consider now the short exact sequence
\begin{displaymath}
\shex{\pi^*\Omega^1_{\P^1}}{\Omega^1_X}{\Omega^1_{X/\P^1}}
\end{displaymath}
and the corresponding long exact sequence of higher direct images under
$\pi$. We have that $\pi_*\Omega^1_{X/\P^1}=R^2\pi_*\Omega^1_{X/\P^1}=0$ as this
is true fibrewise which makes the long exact sequence split up into the
isomorphism $\Omega^1_{\P^1} \riso \pi_*\Omega^1_X$ and the long exact sequence
\begin{displaymath}
0 \to R^1\pi_*\Omega^1_X \mapright{} R^1\pi_*\Omega^1_{X/\P^1} \mapright{}
\Omega^1_{\P^1}\Tensor R^2\pi_*\cO_X \mapright{} R^2\pi_*\Omega_X \to 0.
\end{displaymath}
Now the map $R^1\pi_*\Omega^1_{X/\P^1} \to \Omega^1_{\P^1}\Tensor R^2\pi_*\cO_X$
is the map induced by the Gauss-Manin connection and as K3-surfaces fulfill the
infinitesimal Torelli theorem and if the Kodaira-Spencer map of $\pi$ is
non-zero this map is non-zero. In particular, $R^2\pi_*\Omega_X$ is then a
torsion sheaf. Furthermore, $\Omega^1_{\P^1}\Tensor R^2\pi_*\cO_X \riso
\cO_{\P^1}(-2)\Tensor\cO_{\P^1}(-2)=\cO_{\P^1}(-4)$ so that there is a non-zero
homomorphism $R^1\pi_*\Omega^1_{X/\P^1} \to \cO_{\P^1}(-4)$. Assume for the
moment that the Kodaira-Spencer map is indeed non-zero. It follows
from (\ref{fundamental diagram}) that we have a short exact sequence
\begin{displaymath}
T\Tensor\cO_{\P^1} \to R^1\pi_*\Omega^1_{X/\P^1} \to K \cap F^*K \to 0
\end{displaymath}
and as the composite of $T\Tensor\cO_{\P^1} \to R^1\pi_*\Omega^1_{X/\P^1}$ with
the Kodaira-Spencer map is zero as $\cO_{\P^1}(-4)$ is negative we get an
induced non-zero map $K \cap F^*K \to \cO_{\P^1}(-4)$ which proves
\DHrefpart{KS1}. On the other hand, $K \cap F^*K$ is a subbundle of
$T_0\Tensor\cO_{\P^1}$ so it is a sum of non-positive line bundles which
together gives that $\det(K\cap F^*K) \le -4$ which excludes $p=5$ as in that
case the Kodaira-Spencer map is always non-zero. We have thus proved
\DHrefpart{KS}.

Let $V$ be the subbundle of $K$ generated by its global sections (which is a
trivial subbundle as $K$ is a subbundle of a trivial bundle). Now, as in
particular $V$ is a subbundle of $K$ $F^*V$ is a subbundle of $F^*K$ but as
$F^*K/K \cap F^*K = L^{\tensor -1}\riso \cO_{\P^1}(-2)$ and $F^*V$ is trivial this implies that
$F^*V$ is contained in $K \cap F^*K$ but then it is contained in $V$ and thus
equal to it so that $F^*V=V$ which means that $V=V_0\Tensor\cO_{\P^1}$, where
$V_0 \subset T_0$ is an $\Z/p$-subvector space. This means that $V_0$ is
contained in $K_s$ for every $s \in \P^1(\k)$. For those $s$ for which $T$
equals $\Pic(X_s)$ -- true for all but a finite number of points -- we have that
$T$ equals the Tate module (cf.\ \lcite[Def.\ 3.2]{Og79}) by \lcite[Cor.\
1.6]{Og79} and hence $K_s$ contains $V_0$ which by \lcite[3.12.3]{Og79} shows
that $V_0=0$ and thus that $H^0(\P^1,K)=0$. This implies that $K$ is a sum of
strictly negative line bundles which gives that $\deg\det K \le
-\rk(K)=-\sigma_0$ but $\deg\det K$ equals $-4$ if $p=2$ and $-2$ if $p=3$. We
thus get \DHrefpart{van}.

To prove \DHrefpart{ii} we use that $\pi_*\Omega^1_X = \Omega^1_{\P^1}$ and
that $R^1\pi_*\Omega^1_X$ is the kernel of the Kodaira-Spencer map
$R^1\pi_*\Omega^1_{X/\P^1} \to \cO_{\P^1}(-4)$ so that in particular
$H^0(R^1\pi_*\Omega^1_X)=H^0(R^1\pi_*\Omega^1_{X/\P^1})$. Now, we have
$R^1\pi_*\Omega^1_{X/\P^1}=T'\Tensor\cO_{\P^1}\Dsum H''$ and $H''$, by
(\ref{fundamental diagram}) into an exact sequence
\begin{displaymath}
\shex{T_0\Tensor\cO_{\P^1}/K+F^*K}{H''}{K\cap F^*K}
\end{displaymath}
and as $H^0(\P^1,K\cap F^*K)=0$ we get that
$H^0(R^1\pi_*\Omega^1_{X/\P^1})=T'\Tensor\k\Dsum H^0(\P^1,T_0/K+F^*K)$. However,
as $T_0/K+F^*K$ is a vector bundle on $\P^1$ generated by global sections we have
that
$h^0(\P^1,T_0/K+F^*K)=\deg(T_0/K+F^*K)+\rk(T_0/K+F^*K)=m+\sigma_0-1$. Further,
when the Kodaira-Spencer map is non-zero we get that $R^2\pi_*\Omega_X$ is
torsion. Finally, the Leray spectral sequence for $\Omega^1_X$ thus gives that
$H^0(X,\Omega^1_X)=H^0(\P^1,\pi_*\Omega_X)=0$,
$H^0(X,\Omega^3_X)=H^1(\P^1,R^2\pi_*\Omega_X)=0$ and an exact sequence
\begin{displaymath}
\shex{H^1(\P^1,\pi_*\Omega_X)}{H^1(X,\Omega^1_X)}{H^0(R^1\pi_*\Omega^1_X)}.
\end{displaymath}
Combining these results the statement follows.

As for \DHrefpart{iii} we have a short exact sequence of connections
\begin{displaymath}
\shex{T_0\Tensor\cO_{\P^1}/F^*K}{H'}{F^*K}
\end{displaymath}
giving an exact sequence of sheaves of horisontal sections
\begin{displaymath}
0 \to T_0\Tensor\cO_{\P^1}/K \longrightarrow \cH^0(H',\nabla) \longrightarrow K
\end{displaymath}
and hence an exact sequence of global horisontal sections
\begin{displaymath}
0 \to H^0(\P^1,T_0\Tensor\cO_{\P^1}/K) \longrightarrow H^0(\P^1,(H',\nabla))
\longrightarrow H^0(K),
\end{displaymath}
but $H^(K)=0$ this shows that the global horisontal sections of $H'$ equals
$H^0(\P^1,T_0\Tensor\cO_{\P^1}/K)$ and as $T_0\Tensor\cO_{\P^1}/K)$ is generated
by global sections (and is a sheaf on $\P^1$) we have
\begin{displaymath}
h^0(T_0\Tensor\cO_{\P^1}/K)=\chi(T_0\Tensor\cO_{\P^1}/K)=\deg(T_0\Tensor\cO_{\P^1}/K)+\rk(T_0\Tensor\cO_{\P^1}/K)=\delta+\sigma_0.
\end{displaymath}
The result then follows from (\ref{DR cohomology}).
\end{proof}
\end{proposition}
\begin{lemma}\label{W/p2}
Let $R$ be a complete local topologically finitely generated $\W(\k)$-algebra
such that $R/pR$ is formally $\k$-smooth. Then either $pR=0$ or $R$ has $\W(\k)/p^2$ as a
quotient.
\begin{proof}
We may write $R$ as a quotient of $\pow[\W(\k)]{t_1,\dots,t_n}$ such that the
quotient map induces an isomorphism $\pow[\k]{t_1,\dots,t_n} \riso R/pR$. Hence
all the elements in the kernel $I$ of the quotient map are divisible by $p$. If
$p$ lies in the ideal we are finished so we may assume that every element in $I$
is $p$ times a non-unit. Consider now the $\W(\k)$-homomorphism
$\pow[\W(\k)]{t_1,\dots,t_n} \to \W(\k)$ that takes the $t_i$ to $p$. By
assumption $I$ maps into $p^2\W(\k)$ and hence we get an induced map $R \to
\W(\k)/p^2$.
\end{proof}
\end{lemma}
The following result will be used to compute the miniversal deformation space of
$X$.
\begin{proposition}\label{deformation relation}
The forgetful map from the formal deformation functor $D_\pi$, of $\k$-deformations of
\map{\pi}{X}{\P^1} to the formal deformation functor $D_X$ of $X$ is formally
smooth. More precisely $D_X$ is the quotient of $D_\pi$ by the action of the
formal completion of the completion of $\PGL_2$.
\begin{proof}
Given a deformation $\Pi$ of $\pi$ $S$ and an extension of the induced
deformation $\cX$ of $X$ over a nil-thickening $S'$ of $S$, as $H^2(X,\co_X)=0$ we may
lift $\Pi^*\cO_{\P^1}(1)$ to a line bundle $\cL$ on $\cX$. As
$H^1(X,\pi^*\cO_{\P^1}(1))=0$  we get that $H^0(\cX,\cL)$ is free of rank two and
maps surjectively onto $H^0(\cX_{|S},\cL_{|S})$. Lifting a basis for it gives a
map to $\P^1_{S'}$ lifting $\Pi$. The last part follows from the fact that the
only choices made are exactly that of a basis for $H^0(\cX,\cL)$.
\end{proof}
\end{proposition}
In characteristic $2$ the general result of Deligne-Illusie (cf,
\cite{deligne87::relev+rham}) on the degeneration of the Hodge to de Rham
spectral sequence is not applicable to threefolds. We record some further
consequences of their arguments that will act as substitutes. 
\begin{lemma}\label{de Rham degeneration}
\part For a smooth and proper $\k$-variety $X$, $\k$ a perfect of characteristic
$p>0$, that is liftable to $\W_2(\k)$ the differentials
\map{d_2}{H^i(X,\Omega^1_{X'})}{H^{i+2}(X,\Cal O_{X'})} of the conjugate spectral
sequence are zero, $X'$ being the pullback of $X$ by the Frobenius map of
$\Sp\k$.

\part If furthermore, $h^0(\Omega^2_X)=h^2(\Cal O_X)=0$, then the dimension of
$H^2_{DR}(X/\k)$ equals $h^1(\Omega^1_X)$.
\begin{proof}
The $d_2$-differentials in question are obtained in the following manner: By the
Cartier isomorphism the canonical truncation $\tau_{\le 1}\Omega^\cdot_X$ fits
into the distinguished triangle
\begin{displaymath}
\to \Cal O_{X'} \to \tau_{\le 1}\Omega^\cdot_X \to \Omega^1_{X'}[-1] \to \Cal
O_{X'}[1] \to
\end{displaymath}
and $d_2$ is simply $H^*(X,-)$ applied to the map $\Omega^1_{X'}[-1] \to \Cal
O_{X'}[1]$. However, by \cite[Thm.\ 2.1]{deligne87::relev+rham} this map is
zero. The second part now follows directly from the first.
\end{proof}
\end{lemma}
\begin{subsection}{The case $\sigma_0(\pi)=2$}

We shall now consider the case of characteristic $\sigma_0=2$. Note that when
$\Char\k=3$ we may, as in \cite{schroeer03::some+calab+yau+witt}, construct, for
each polarisation $\cL$ on the product of supersingular elliptic curves with
$K(\cL)\riso \alpha_3\times\alpha_3$, a family of Kummer surfaces over $\P^1$
associated to the Moret-Bailly family associated to $\cL$. For the sake of
convenience we shall call such a family a \Definition{Moret-Bailly family}. It
follows from \cite[Thm.\ 7.15]{Og79} that the moduli space of supersingular
K3-surfaces marked with a K3-lattice with $\sigma_0 \le 2$ is the union of these
families. When $\Char\k=2$ we can instead follow
\cite{schroeer03::some+calab+yau+witt} and use the generalised Kummer
construction. As we have a Torelli theorem for supersingular K3-surfaces in
characteristic $2$ (cf.\ \cite{rudakov81::surfac+k3} and \cite[p.\
388]{ogus83::ae+torel+k3}) and the coarse moduli space of $\sigma_0\le 2$
K3-periods is irreducible (cf.\ \cite[Thm.\ 3.21]{Og79}) we get that the
irreducible components of the moduli space of $\sigma_0 \le 2$ K3-surfaces are
all obtained in this way. (This can presumably be proved by a variant of Ogus'
proof in odd characteristic.) We shall again call such a family a
\Definition{Moret-Bailly} family.
\begin{theorem}\label{char=3}
Let \map{\pi}{X}{\P^1} be a smooth and proper map over an algebraically closed
field of characteristic $3$, the fibres of which are supersingular K3-surfaces
and for which $\omega_X$ is trivial. 

\part[i] $\pi$ is a Moret-Bailly family, in particular the Artin invariant of its
generic fibre is $2$.

\part[ii] The Hodge numbers of $X$ are as follows
\begin{displaymath}
\begin{array}{|c||c|c|c|c|}
\hline
3&1&0 &0 &1\\
2&0&0 &24&0\\
1&0&24&0 &0\\
0&1&0 &0 &1\\
\hline\hline
ji&0&1&2&3\\
\hline
\end{array}
\end{displaymath}

\part[iii] The $\k$-dimension of $H^i_{DR}(X/\k)$ equals $b_i$.

\part Any $X$ is arithmetically rigid, i.e., its $W(\k)$-formal deformation
space equals $\Sp\k$.
\begin{proof}
We notice to begin with that if the Kodaira-Spencer is zero then $\pi$ (with a
given a marking) is the pullback over the Frobenius map of another fibration
which would force $R^2\pi_*\cO_X$ to also be the pullback along the Frobenius
map. This, however, is not possible as the degree of $R^2\pi_*\cO_X$, $-2$, is
not divisible by $3$.

By Proposition \ref{quotient possibilities} we see that the generic Artin
invariant $\sigma_0(\pi)$ must be equal to $2$ and hence $\pi$ marked by a
marking of the Picard lattice of its generic fibre is the pullback of a
Moret-Bailly family by a map \map{f}{\P^1}{\P^1}. Now, by \cite[Prop.\
4.2]{schroeer03::some+calab+yau+witt} $\omega_{X'}$ for a Moret-Bailly family
\map{\pi '}{X'}{\P^1} is trivial and hence $\omega_{X'/\P^1}\riso \pi
'^{*}\cO(2)$ and as this is also true for $\omega_{X/\P^1}$ we conclude that $f$
has degree $1$ which gives \DHrefpart{i}.

We get the Hodge numbers from Propositions \ref{prel cohomology} and
\ref{quotient possibilities}.

As for \DHrefpart{iii} we get the result for $i=0,1,5,6$ from \DHrefpart{ii}
and from duality and the fact that the Euler characteristic is $48$ it is enough
to determine $\dim_{\k}H^2_{DR}(X/\k)$ which we get from (\ref{DR cohomology1}).

Finally, from Lemma \ref{W/p2} we get that if $X$ is not arithmetically rigid
then there is a non-trivial deformation of $X$ over $\W(\k)/9$ but this
contradicts \cite[Cor.\ 2.4]{deligne87::relev+rham} as the Hodge to de Rham
spectral sequence of $X$ does not degenerate.
\end{proof}
\end{theorem}
We now turn to the case of characteristic $2$.
\begin{theorem}\label{sigma=2}
Let \map{\pi}{X}{\P^1} be a smooth and proper map over an algebraically closed
field of characteristic $2$, the fibres of which are supersingular K3-surfaces
with the generic fibre of Artin invariant $2$ and for which $\omega_X$ is
trivial.

\part $\pi$ is the pullback of a Moret-Bailly family by a map $\P^1 \to \P^1$,
called the \Definition{characteristic map} of
degree $2$ and any such pullback is of the desired form. Hence the
Kodaira-Spencer map is zero precisely when the characteristic map is inseparable
(i.e., isomorphic to the Frobenius map).

\part[ii] When the Kodaira-Spencer map is non-zero the Hodge numbers of $X$ are
as follows
\begin{displaymath}
\begin{array}{|c||c|c|c|c|}
\hline
3&1&0 &0 &1\\
2&0&2 &26&0\\
1&0&26&2 &0\\
0&1&0 &0 &1\\
\hline\hline
ji&0&1&2&3\\
\hline
\end{array}
\end{displaymath}
and when it is zero they are as follows
\begin{displaymath}
\begin{array}{|c||c|c|c|c|}
\hline
3&1&3 &0 &1\\
2&0&5 &26&0\\
1&0&26&5 &0\\
0&1&0 &3 &1\\
\hline\hline
ji&0&1&2&3\\
\hline
\end{array}
\end{displaymath}

\part[iii] The $\k$-dimensions of $H^i_{DR}(X/\k)$ equals $1$, $0$, $25$, $4$,
$25$, $0$, and $1$ for $i=0,1,2,3,4,5,6$.

\part The $W(\k)$-formal deformation
space of $X$ is smooth over $\k$ of dimension $2$ if the Kodaira-Spencer of the
pencil is non-zero smooth over $\k$ of dimension $5$ if it is.
\begin{proof}
The map from the base of $\pi$, $\P^1$, is irreducible and hence it maps to one
component of the moduli space of marked K3-surfaces with Artin invariant $\le
2$, which as we have seen is a Moret-Bailly family. Now, for a Moret-Bailly
\map{\pi '}{X'}{\P^1} we have that $R^1\pi '_*\cO_{X'}\riso
\cO_{\P^1}(-1)$. This can be done directly but the quickest way is to note that
by \cite[Prop. 4.2]{schroeer03::some+calab+yau+witt} the pullback by a map of
degree $2$ of $\P^1$ is $\cO_{\P^1}(-2)$. Hence, the characteristic map is of
degree $2$ and conversely any map $\P^1 \to \P^1$ has a pullback the total space
of which is Calabi-Yau. Clearly the Kodaira-Spencer map is zero precisely when
the characteristic map is inseparable as that of a Moret-Bailly family as
non-zero. However, a map of degree $2$ is inseparable precisely when it is the
Frobenius map.

The second part follows when the Kodaira-Spencer map is non-zero just as in the
characteristic $3$ case from Propositions \ref{prel cohomology} and
\ref{quotient possibilities}. When the Kodaira-Spencer map is zero we get from
the long exact sequence of cohomology that $R^1\pi_*\Omega^1_X \riso
R^1\pi_*\Omega^1_{X/\P^1}$ and $\cO_{\P^1}(-4) \riso R^2\pi_*\Omega^1_X$. This and
the fact that $R^1\pi_*\Omega^1_{X/\P^1}$ splits as an orthogonal direct sum
$T'\Tensor\cO_{\P^1}\perp T'^\perp$ and that the latter factor fits into an exact sequence
\begin{displaymath}
\shex{\cO_{\P^1}(6)}{T'^\perp}{\cO_{\P^1}(-6)}
\end{displaymath}
give the Hodge numbers.

As for the de Rham-cohomology the computation of $\dim_\k H^2_{DR}(X/\k)$
follows from (\ref{DR cohomology1}) and the rest then follows from duality.

Finally, the space of maps $\P^1 \to \P^1$ of degree $2$ is smooth (of dimension
$3$) and hence it follows from Proposition \ref{deformation relation} that a
miniversal formal $\k$-deformation space is smooth. On the other hand, as
$T^1_X\riso \Omega^2_X$ we get $H^1(X,T^1_X)$ has dimension $2$ or $5$. In the
case of non-zero Kodaira-Spencer map assume that a miniversal
$\W(\k)$-deformation space does not equal the $\k$-deformation space then by
Lemma \ref{W/p2} there would be a lifting of $X$ to $\W/p^2$.  That however
contradicts that the formulas for the Hodge numbers and dimensions of de
Rham-cohomology shows that the conditions of the second part of Lemma \ref{de
Rham degeneration} are fulfilled yet the conclusion is false.

The case when the Kodaira-Spencer map is zero requires some elaboration. We
start by noticing that if $X$ is liftable to $\W/4$ then so is the pencil $X
\to \P^1$ as first the line bundle $\cO(1)$ giving the map is liftable as
$H^2(X,\cO_X)=0$ and then the sections are liftable as $H^1(X,\cO(1))=0$. Now,
the obstruction for lifting the map (mapping to a fixed lifting of $\P^1$) is an
element of $H^2(X,T_{X/\P^1})$ and we shall show that it is non-zero. Consider
instead the Moret-Bailly family \map{\pi '}{X'}{\P^1} for which we also have an
obstruction in $H^2(X',T_{X'/\P^1})$ to the liftability of the pencil. Now,
$H^2(T_{X'/\P^1})=H^1(\P^1,R^1\pi 'T_{X'/\P^1})$ and the same is true for any
pencil obtained from the Moret-Bailly pencil by taking the pullback along any
non-constant map $\P^1 \to \P^1$. We have just seen that the pullback along a
separable map of degree $2$ is not liftable and hence the obstruction in
$H^2(T_{X'/\P^1})=H^1(\P^1,R^1\pi 'T_{X'/\P^1})$ is non-zero. However, for all
vector bundles $\Cal E$ on $\P^1$ and all non-constant maps \map{f}{\P^1}{\P^1},
the pullback map \map{f^*}{H^1(\P^1,\Cal E)}{H^1(\P^1,f^*\Cal E)} is
injective. For this one writes $\Cal E$ as a sum of line bundles and then only
the case $\cO(n)$ for $n\le -2$ is non-trivial. In that case the cokernel of
$\cO(n) \to f_*f^*\cO(n)$ has just zero as global section.
\end{proof}
\end{theorem}
\end{subsection}
\begin{subsection}{The case $\sigma_0(\pi)=3$}

We are now going to consider the case of $\sigma_0=3$ and for that we shall need
some preliminaries. We are going to identify the $\sigma_0 \le 3$ period space,
which by its very definition is a Deligne-Lusztig type variety for the
orthogonal group of a non-split form on $\F_p^6$ with the degree $p+1$ Fermat
surface in $\P^3$ which is a Deligne-Lusztig type variety for the unitary group
$\symb{SU}_3$ over $\F_p$. That there is such an identification follows from the
coincidence $A_3=D_3$. However, it is probably easier to construct this
correspondence directly and that concrete construction also simplifies some
identifications we shall have to make.

To begin with we need to explicitly handle the Pl\"ucker embedding in
characteristic $2$ which means that we need not just a scalar product but also a
quadratic form on $\Lambda^2\k^4$. This can be done directly (cf.\ \cite[Lemma
8.3]{eisenbud00::enriq+pfaff}) but we find the following construction more
natural. For that we extend the notion of divided powers to graded ideals in
$\Z/2$-graded commutative algebras. The extension is straightforward; we demand
that all divided powers of degree $\ge 2$ be zero and also the usual conditions
with appropriately inserted signs. In particular the even part of such an ideal
considered as an ideal in the even part of the algebra will have a divided power
structure in the usual sense.
\begin{proposition}
If $M$ is a module over a commutative ring $R$ then 
$\Lambda^+M$ has a unique divided power structure.
\begin{proof}
A divided power structure is unique as $\Lambda^*M$ is generated by odd
elements. To prove existence we first write $M$ as a quotient of a free module
$F$ and then $\Lambda^*M$ is the quotient of $\Lambda^*F$ by the ideal generated
by the kernel of the quotient map and it is clear that a divided power structure
on $\Lambda^*F$ will preserve that ideal. Hence we may assume that $M$ is free
and then by a limit argument to it being also finitely generated and then one
reduces to the module having rank $1$ by using that the graded tensor product of
divided powers has a compatible divided power. The case of rank $1$ is trivial.
\end{proof}
\end{proposition}
In particular we get the promised quadratic form on $\Lambda^2\k^4$ with values
in $\Lambda^4\k^4$ and it is clear that $\Lambda^2U$ of a $2$-dimensional
subspace $U$ of $\k^4$ is isotropic with respect to it. We shall however need a
twisted version of this construction. Hence let $V$ be a $4$-dimensional vector
space over $\F_{p^2}$, $p$ a prime, and $\skpr--$ a non-degenerate hermitian
form on $V$. It induces a hermitian form on $\Lambda^4V$ and there is an element
$w \in \Lambda^4V$ such that $\skpr{w}{w}=1$ and we fix one such $w$.We now
define a semi-linear automorphism \map{\overline{(-)}}{\Lambda^2V}{\Lambda^2V} by $u\wedge v =
\skpr{\overline u}{v}w$, where $\skpr--$ is the hermitian scalar product on
$\Lambda^2V$ induced by that on $V$. We claim that this is an involution. Indeed, this is
independent of the choice of $w$ as any other such choice is of the form
$\lambda w$ with $\lambda\overline{\lambda}=1$ and for such a different choice
the semi-linear map is $u \mapsto \lambda \overline{u}$ but
$\lambda\overline{\lambda\overline{u}}=\lambda\overline{\lambda}\overline{\overline{u}}=\overline{\overline{u}}$.
We may hence choose an orthonormal basis $e_1,e_2,e_3,e_4$ for $V$ and then
$\overline{e_i\wedge e_j}$ equals $\pm e_k\wedge e_\ell$, where $\{k,\ell\}$ is
the complement of $\{i,j\}$ and the sign is chosen so that $e_i\wedge e_j\wedge
e_k\wedge e_\ell = \pm e_1\wedge e_2\wedge e_3\wedge e_4$. This shows that we
have an involution.

We may now define an $\F_p$-vector space $W$ as the fixed points of the
involution. We then get a quadratic $\F_p$-valued form $\psi$ on $W$ by
$\gamma_2(u)=\psi(u)w$. (This construction is essentially independent of the
choice of $w$ as if we choose $\lambda w$ instead we may find $\mu \in \F_{p^2}$
such that $\lambda=\mu^{p-1}$ and then multiplication by $\mu$ takes the fixed
points for one involution to the fixed points for the other.) It follows from
general principles that this quadratic form is the non-split form but it can
also easily be verified directly by choosing an explicit basis for $W$. We do it
only in the case that interests us, namely when $p=2$. In that case we choose
$\zeta \in \F_4\setminus \F_2$ and then $\zeta^3=1$. Then choosing an
orthonormal basis as before we have that $e_i\wedge e_j+e_k\wedge e_\ell$ and
$\zeta e_i\wedge e_j+\zeta^{-1}e_k\wedge e_\ell$, for $\{i,j\}$ running over
$\{1,2\}$, $\{1,3\}$, and $\{1,4\}$ and $k$ and $\ell$ as before, form a basis
for $W$. Then it is clear that the spaces spanned by $e_i\wedge e_j+e_k\wedge
e_\ell$ and $\zeta e_i\wedge e_j+\zeta^{-1}e_k\wedge e_\ell$ for the different
pairs $\{i,j\}$ are orthogonal to each other and each of them is a non-split
plane. The orthogonal direct sum of three non-split planes is non-split however.

Assume now that $S$ is an $\F_{p^2}$-scheme and $\cL \subset \cO_S^4$ is a
sub-line bundle. Then $\cL\wedge\cO_S^4/\cL \subset \Lambda^2\cO_S^4$ is a
maximal isotropic sub-vector bundle (with respect to the quadratic form $\psi$
which of course is isomorphic to the standard form as $S$ is an
$\F_{p^2}$-scheme). Similarly, if $\cH \subset \cO_S^4$ is a rank $3$
sub-vector bundle then $\Lambda^2\cH \subset \Lambda^2\cO_S^4$ is again a maximal
isotropic subbundle. In this way we get an isomorphism between $\P^3_{\F_{p^2}}$
and $\check{\P}^3_{\F_{p^2}}$ and the two irreducible components over $\F_{p^2}$
of maximal isotropic subspaces of $(W,\psi)$. We shall however need to identify
the action of the Galois group of $\F_{p^2}/\F_p$ on these components (which of
course is trivial if we instead had considered a split form).
\begin{proposition}\label{Galois action}
\part[i] Let $S$ be an $\F_{p^2}$-scheme, $V = \F_{p^2}^4$, $0 \ne \Lambda^4V$,  and $\cL \subset
V\Tensor \cO_S$ a sub-line bundle. Let $\cH \subset V^*\Tensor \cO_S$ be the
annihilator of $\cL$. We may identify $\Lambda^2V$ and $\Lambda^2V^*$ using first
the canonical identification $\Lambda^2V^*=(\Lambda^2V)^*$ and then the
identification $(\Lambda^2V)^*=\Lambda^2V$ given by $u(v)=u\wedge v/w$. Under
this identification $\cL\wedge V\Tensor \cO_S/\cL \subset \Lambda^2V\Tensor\cO_S$
and $\Lambda^2\cH \subset \Lambda^2V^*\Tensor\cO_S$ correspond to each other.

\part Let $V$ be a $4$-dimensional vector space over $\F_{p^2}$ with a
non-degenerate hermitian form and let $W$ be the $\F_p$-form of $\Lambda^2V$
given by the hermitian form and an element of $\Lambda^4V$ of square $1$ and let
$\psi$ be the corresponding quadratic form on $W$. Let $M$ be the moduli space of
maximal isotropic subspaces of $W$. Then $M$ scalar extended to $\F_{p^2}$ is the
disjoint union of $\P(V)$ and $\check{\P}(V)=\P(V^*)$ with the universal maximal
isotropic subspace on $\P(V)$ given by $\cL\wedge V \subset \Lambda^2V$ on
$\P(V)$ and the one on $\check{\P}(V)$ by $\Lambda^2\cH \subset \Lambda^2V$, where $\cL \subset
V\Tensor\cO_{\P(V)}$ and $\cH \subset V\Tensor\cO_{\check{\P}(V)}$ are the
tautological subbundles. The action of $\sigma$, the $p$'th power map $\F_{p^2}
\to \F_{p^2}$, on $M$ takes $\cL$ to the annihilator of $\sigma^*\cL \subset
\sigma^*V$ under the hermitian form.
\begin{proof}
For the first part we use that $\Lambda^2\cH \subset \Lambda^2V^*\Tensor\cO_S$
is the kernel of the contraction map
\map{\ell\contract}{\Lambda^2V^*\Tensor\cO_S}. Now, the identification
$\Lambda^2V^* \to (\Lambda^2V)^*$ is given by $f\wedge g(u\wedge
v)=f(u)g(v)-f(v)g(u)$ and the contraction is a derivation so that
$\ell\contract(f\wedge g)=f(\ell)g-g(\ell)f$ so that $(\ell\contract f\wedge
g)(u)=f(\ell)g(u)-g(\ell)f(u)=f\wedge g(\ell\wedge u)$, which means that the
contraction transfers to $\ell\contract F(U):= F(\ell\wedge u)$ for $F \in
(\Lambda^2V)^*$. Now \DHrefpart{i} means that the image of the composite
\begin{displaymath}
\cL\wedge V\Tensor \cO_S/\cL \subset \Lambda^2V\Tensor\cO_S \riso (\Lambda^2V)^*\Tensor\cO_S
\end{displaymath}
equals the image of the composite
\begin{displaymath}
\Lambda^2\cH \subset \Lambda^2V^*\Tensor\cO_S \riso (\Lambda^2V)^*\Tensor\cO_S
\end{displaymath}
but as both are subbundles of equal rank it is enough to show that the image of
the first map lies in the image of the second. This follows from the above
formula as if $F \in (\Lambda^2V)^*$ is defined by $F(v) := \ell\wedge u\wedge
v/w$ then $(\ell\contract F)(v)= F(\ell\wedge v) = (\ell\wedge\ell\wedge u\wedge
v)/w = 0$.

Turning to the second part, by the definition of the involution on $\Lambda^2V$
it follows directly from the first part.
\end{proof}
\end{proposition}
With these preliminaries we can now identify the period space.
\begin{proposition}\label{period space}
Consider the hermitian form on $\F_{p^2}^4$ for which the standard basis is
orthonormal and let $X_p \subset \P^3_{\F_{p^2}}$ be the hypersurface defined by
$x^{p+1}+y^{p+1}+z^{p+1}+w^{p+1}=0$. Let $\psi$ be the quadratic form on the
$\F_p$-form $W$ of $\Lambda^2\F_{p^2}^4$ defined by the hermitian form and
$e_1\wedge e_2\wedge e_3\wedge e_4 \in \Lambda^4\F_{p^2}^4$. 

\part The scalar extension to $\F_{p^2}$ of the period space of maximal
isotropic subbundles $K$ for $\psi$ for which $K\cap F^*K$ is a subbundle of
rank $3$ can be identified with two copies of $X_p$; the universal subbundle on
one copy being $\cL\wedge \cO_{X_p}/\cL \subset
\Lambda^2\F_{p^2}^4\Tensor\cO_{X_p}$ and $\Lambda^2\cH \subset
\Lambda^2\F_{p^2}^4\Tensor\cO_{X_p}$ on the other, where $\cL \subset
\cO_{X_p}^4$ is the (restriction of the) universal subbundle and $\cH$ its
annihilator (with respect to the standard scalar product).

\part\label{rational lines} For each $2$-dimensional $\F_{p^2}$-vector space $V$
of $\F_{p^2}^4$ isotropic with respect to the hermitian form, the line $\P^1(V)
\subset \P^3$ will lie in $X_p$ and represents the space of totally isotropic
subspaces in the first component which all contain $\Lambda^2V$. In this way
one gets all the maximal isotropic subspaces of $X_p$ containing a fixed
$1$-dimensional $\F_p$-rational isotropic subspace. We shall call such lines
\Definition{rational lines}. The intersections of two rational lines gives all
the maximal isotropic subspaces of $X_p$ containing a fixed $2$-dimensional
$\F_p$-rational isotropic subspace.

\part\label{tautological determinant} The determinant of the tautological
totally isotropic subbundle in the first component equals $\cL^{\tensor2}$.
\begin{proof}
For the first part we note that if $S'$ is an $\F_p$-scheme and $S$ is its
scalar extension to $\F_{p^2}$, then the (absolute) Frobenius map on $S$ factors
as the composite of the Frobenius map on $\Sp\F_{p^2}$ and the relative
Frobenius for $S \to \Sp\F_{p^2}$. Combined with Proposition \ref{Galois action}
this immediately gives the first part.

As for the second part an $\F_p$-rational $1$-dimensional isotropic subspace $U$
of $\Lambda^2\F_{p^2}^4$ corresponds by the Pl\"ucker embedding to a
$2$-dimensional subspace $V$ of $\F_{p^2}^4$ and by Proposition \ref{Galois action}
that $U$ is $\F_p$-rational is equivalent to $V$ being isotropic.

Finally, we have that the universal subbundle is isomorphic to
$\cL\Tensor\F_{p^2}^4/\cL$ and hence has determinant
$\cL^{\tensor3}\Tensor\cL^{\tensor-1}$.
\end{proof}
\end{proposition}
\begin{remark}
The same proof as in the proof of $H^0(\P^1,K)=0$ in Proposition \ref{quotient
possibilities} shows that for any family $X \to S$ of supersingular K3-surfaces
$H^0(S,K)$ is $\F_p$-rational. This shows in particular that \emph{all} lines on
$X_p$ are constructed in the way of (\ref{rational lines}).
\end{remark}
We are now ready to give the description of the pencils with $\sigma_0=3$. We
note that for $p=2$ the period space is, according to Proposition \ref{period
space}, the Fermat cubic $X_2$. If $\ell$ is one of the $27$ lines on $X_2$ and
$H$ is a hyperplane containing $\ell$, then $H$ intersects $X_2$ in $\ell$ and a
residual quadric $Q_H$. Except for the $5$ cases, for each line $\ell$, when
$Q_H$ is a union of two lines, $Q_H$ is a smooth plane quadric and if
\map{\pi}{X}{Q_H} is the pullback of the universal family of K3-surfaces on
$X_2$ it follows from (\ref{line equation}) and (\ref{tautological determinant})
that $R^2\pi_*\cO_X \riso \cO_{Q_H}(-2)$ so that the total space of $X$ is a
Calabi-Yau threefold. As $Q_H$ does not lie in the union of the rational lines
the general point will, by (\ref{rational lines}), have $\sigma_0=3$. We shall
call the family given by $Q_H$ the \Definition{residual family} with respect to
$H$.
\begin{theorem}\label{sigma=3}
Let \map{\pi}{X}{\P^1} be a smooth and proper map over an algebraically closed
field of characteristic $2$, the fibres of which are supersingular K3-surfaces
with the generic fibre of Artin invariant $3$ and for which $\omega_X$ is
trivial.

\part The period mapping of $\pi$ maps $\P^1$ isomorphically to a residual
quadric for one of the lines of the Fermat cubic.

\part[ii] The Hodge numbers of $X$ are as follows
\begin{displaymath}
\begin{array}{|c||c|c|c|c|}
\hline
3&1&0 &0 &1\\
2&0&1 &25&0\\
1&0&25&1 &0\\
0&1&0 &0 &1\\
\hline\hline
ji&0&1&2&3\\
\hline
\end{array}
\end{displaymath}

\part[iii] The $\k$-dimensions of $H^i_{DR}(X/\k)$ equals $1$, $0$, $24$, $2$,
$24$, $0$, and $1$ for $i=0,1,2,3,4,5,6$.

\part[def] The $W(\k)$-formal deformation space of $X$ is smooth over $\k$ of
dimension $1$.
\begin{proof}
The period map maps $\P^1$ into the period space and by irreducibility it lands
in one component and by duality or conjugation we can assume that it lands in
the one given by $\cL\Tensor\cO^4/\cL$. By (\ref{line equation}) and
(\ref{tautological determinant}) the obtained map $\P^1 \to X_2 \subset \P^3$ is
of degree $2$. If the map is not birational onto its image, the image is a
line. However, there are $27$ rational lines (as for instance there are exactly
$27$ isotropic non-zero elements for the non-split form on $\F^6_2$) and hence
all lines in $X_3$ are rational and $\sigma_0(\pi) \le 2$ and not equal to $3$
as is assumed. Hence the image is a curve $Q$ of degree $2$ and thus lies in a
hyperplane $H$. The curve residual to $Q$ in $X_2 \cap H$ is a line and is hence
one of the rational lines so that $\pi$ is a residual family.

The formulas for the Hodge and de Rham numbers now follows from Proposition
\ref{prel cohomology}.

As for \DHrefpart{def}, by varying the hyperplane, each residual family becomes
part of a $1$-dimensional effective family and as $h^1(X,T^1_X)=h^{21}=1$ we see
that the $\k$-deformation space is smooth of dimension $1$ and then we get that
it equals the $\W(\k)$-deformation space by Proposition \ref{W/p2}, Lemma
\ref{de Rham degeneration}, and the computation of de Rham-cohomology dimensions
and Hodge numbers.
\end{proof}
\end{theorem}
Note that as noticed in both the proposition and the previous remarks all the
lines on $X_3$ come from rational lines.  This means that a general residual
quadric is a deformation of two non-skew rational lines as such a pair is a
residual quadric for a hyperplane through any line on a cubic. Each line gives a
$\sigma_0=2$ family and the fact that the lines intersect means that they have a
marked K3-surface in common (which has $\sigma_0=1$). Hence the non-singular
residual quadrics are deformations of a family over two crossing $\P^1$'s with a
point in common. Hence an alternative way of showing that $\sigma_0=3$-families
exist -- which does not make use of period theory -- would be by using deformation
theory to show that the family over two crossing $\P^1$'s deform. We shall now
proceed by sketching such an approach. Let $N$ be a K3-lattice with $\sigma_0=3$
and let $V \subset N^*/N$ be a totally isotropic $2$-dimensional subspace and
choose $e_1\ne e_2 \in V \setminus \{0\}$. The inverse images in $N^*$ of
$\Z/2e_i$ give K3-lattices $N_1$ and $N_2$ with $\sigma_0=2$ and $N_1+N_2$ is a
K3-lattice with $\sigma_0=1$. Now, there is an automorphism of $N^*/N$ permuting
$e_1$ and $e_2$ and it can be lifted (cf.\ \cite[Thm.\ 1.14.2]{Ni80}) to an automorphism
$\varphi$ of $N$ that then permutes $N_1$ and $N_2$. Pick now a smooth pencil of
\map{\pi}{X}{\P^1} of supersingular K3-surfaces with $\sigma_0=2$ marked by $N$
in such a way that the image of $N_1$ is the Picard group of a general
fibre. Choose another marking of the pencil given by the chosen one composed
with $\varphi$. Pick a point of $\P^1$ for which the Picard group of the fibre
is given by $N_1+N_2$ and glue it with itself on another copy of the pencil to
get a nodal genus $0$ curve $P$. The family of K3-surfaces as well as the
marking given by the two constructed markings on the two components also
glues to give a family of $N$-marked K3-surfaces over $P$. This then gives a map
\map{f}{C}{\cM^N} from $C$ to the moduli space of $N$-marked K3-surfaces. This
map is an immersion. This is clear outside of the singular point of $C$ and is
true at the singular point because in one component the marking has Picard group
equal to $N_1$ and on the other equal to $N_2$. Thus we have an exact sequence
\begin{displaymath}
\shex{L}{f^*\Omega_{\cM^N}^1}{\Omega_C^1}
\end{displaymath}
so that the cotangent complex $T^1_f$ of $f$ has the form $L[-1]$. Furthermore,
$\Omega_{\cM_N}$ is a vector bundle and $\Omega_C^1$ has projective dimension
$1$ so that $L$ is a vector bundle and having rank $1$ it is a line
bundle. Hence we can compute by taking determinants of the exact
sequence. Furthermore, a line bundle on $C$ is determined by its degree on each
of the two components of $C$. To begin with, the degree of $\Omega_C^1$ is most
easily computed by realising $C$ as two distinct lines in the plane and from
that it is clear that its degree on each component is $-1$. On the other hand,
from \cite[Cor. 5.4]{Og79} it follows that $\Omega_{\cM_N}^1$ pulled back to either of the
components of $C$ is $\Hom(K\cap F^* K,F^*K/K\cap F^* K)$ and as we have seen
$F^*K/K\cap F^* K \riso \cO_{\P^1}(-1)$ and $K\cap F^* K$ is an extension of
$\cO_{\P^1}\Tensor e_i$ by the similar $K\cap F^* K$ for the
$\sigma_0=2$-marking which again has been seen to be $\cO_{\P^1}(-1)$. We
therefore get that $L \riso \cO_C$. Now first order deformations of the map $f$,
with fixed target space, are given by elements of $H^1(C,T^1_f)\riso
H^0(C,\cO_C)\riso \k$ and obstructions to deformations are given by
$H^2(C,T^1_f)\riso H^1(C,\cO_C)=0$. This shows that the deformation space of $f$
is smooth of dimension $1$ but it is easily seen that deformations that keep the
singularity over a smooth base are trivial (it is not difficult to directly show
that also first order deformations are).
\end{subsection}
\begin{subsection}{The case $\sigma_0(\pi)=4$}

The possibility of there being a pencil with $\sigma_0=4$ is easily taken care of.
\begin{theorem}\label{sigma=4}
There is no smooth and proper map over an algebraically closed field of
characteristic $2$, the fibres of which are supersingular K3-surfaces with the
generic fibre of Artin invariant $4$ and for which $\omega_X$ is trivial.
\begin{proof}
Indeed, we have that $\deg\det K=-4$ and $K$ is strictly negative of rank $4$
which forces $K \riso \cO_{\P^1}(-1)^4$. This gives $F^*K \riso
\cO_{\P^1}(-2)^4$ and as $F^*K/K \cap F^*K \riso \cO_{\P^1}(-2)$ we get that $K
\cap F^*K$ is isomorphic to $\cO_{\P^1}(-2)^3$. This however means that there is
no non-zero map $K \cap F^*K \to \cO_{\P^1}(-4)$ and hence by Proposition
\ref{quotient possibilities} the Kodaira-Spencer map is inseparable so that $\pi$
is the pullback by the Frobenius map of a family \map{\pi '}{X'}{\P^1}. For this
map we will have $K = F^*K'$ which is impossible as $K \riso \cO_{\P^1}(-1)^4$.
\end{proof}
\end{theorem}
\end{subsection}
\begin{proofof}{Theorem \ref{Pencil}}
The theorem follows from Theorems \ref{char=3}, \ref{sigma=2}, and \ref{sigma=3}.
\end{proofof}
\begin{subsection}{Crystalline and de Rham-Witt cohomology}

We finish by computing the crystalline cohomology of the K3-pencils. As for the
Hirokado example we haven't been able to completely determine the de Rham-Witt
cohomology, contrary to that case we are able to completely determine
crystalline cohomology and from that it is possible to determine the de
Rham-Witt cohomology up to finite ambiguity. We have however not been able to
completely determine this cohomology and we shall restrict ourselves to the case
where we can avoid ambiguity.
\begin{proposition}
Let \map{\pi}{X}{\P^1} be a smooth map over a field $\k$ of characteristic $p$
the fibres of which are supersingular K3-surfaces with $\sigma_0$ the Artin
invariant of the general fibre and whose total space is a Calabi-Yau threefold.

\part We have $H^0(X/\W)=H^6(X/\W)=\W$, $H^1(X/\W)=H^5(X/\W)=0$,
$H^2(X/\W)\riso\W^{23}$, $H^3(X/\W)\riso(\W/p\W)^k$, and
$H^3(X/\W)\riso(\W/p\W)^k\Dsum\W^{23}$, where $k=0$ if $p=3$ and $k=4-\sigma_0$
if $p=2$.

\part $R\Gamma(X,\W\Omega^\cdot_X)$
as a direct sum 
\begin{displaymath}
\W\Dsum\W^{23}[-1](-1)\Dsum M[-3] \Dsum\W^{23}[-2](-2)\Dsum\W[-4](-4)
\end{displaymath}
and $M$ is a torsion complex with the properties:
\begin{enumerate}
\item The Hodge numbers are $h^{0,0}(M)=h^{3,-3}(M)=1$,
$h^{-2,-1}(M)=h^{-1,-2}(M)=h^{1,1}(X)-23$, $h^{1,0}=h^{2,-3}=h^{1,3}(X)$ and
$h^{-1,-1}(M)=h^{-2,-2}(M)=h^{1,2}(X)$ and the rest are zero.

\item $T^{0,0}(M)=T^{1,-1}=1$ and the rest are zero.

\item We have $D(M) = M[-3](3)$.
\end{enumerate}

\part When $p=3$ we have that $M$ is isomorphic to the $N_1$ of Lemma \ref{Nk}.
\begin{proof}
The first part follows from the computation that has been done of the de
Rham-cohomology once we know that the torsion is killed by $p$. For that we
consider the Leray spectral sequence for $\pi$ and crystalline cohomology. We
have that $R^i\pi_*\cO_{X/W}$ is zero if $i \ne 0,2,4$, is isomorphic to
$\cO_{X/W}$ if $i=0,4$. For $i=2$ it is the K3-crystal (cf.\ \scite[\S5]{Og82})
associated to $\pi$. As such it contains a Tate crystal (cf.\ \lcite[Prop.\
5.2]{Og82}) $E$ with $E \subset R^2\pi_*\cO_{X/\W} \subset E^*$ and in
particular $R^2\pi_*\cO_{X/\W}/E$ is killed by $p$ (as for instance $E^*/E$
is). Now, from the short exact sequence of crystalline sheaves
\begin{displaymath}
\shex{E}{R^2\pi_*\cO_{X/\W}}{R^2\pi_*\cO_{X/\W}/E}
\end{displaymath}
we get a long exact sequence in cohomology. Furthermore, as
$\Gamma(\P^1,\underline{\Pic}(X/\P^1))\Tensor \cO_{\P^1/W} \to E$ is an
isomorphism, $E$ is isomorphic to a sum of $22$ copies of $\cO_{\P^1/W}$ and so
in particular $H^1(\P^1,E)=0$ by the computation of crystalline cohomology of
the projective line. From this we conclude that $H^1(\P^1,R^2\pi_*\cO_{X/\W})$ injects into
$H^1(\P^1,R^2\pi_*\cO_{X/\W}/E)$ and is hence killed by $p$. By the Leray
spectral sequence we have $H^3(X/\W)=H^1(\P^1,R^2\pi_*\cO_{X/\W})$.

The proof that we may split up $R\Gamma(X,\W\Omega^\cdot_X)$ is the same as the
proof of the same statement in Theorem \ref{dRW}. The Hodge numbers can then be
read off directly from the known Hodge numbers for the Calabi-Yau threefolds,
the formula for the $T$-numbers is obtained as in (\ref{dRW}) and the duality
formula follows from duality for $R\Gamma(X,\W\Omega^\cdot_X)$.

Finally, in the case when $p=3$ we get that $M$ is $\cs$-acyclic and hence is
without finite torsion. We have in particular that $H^0(M)$ is a $1$-dimensional
domino and hence is isomorphic to $\domU_k$ for some $k$. As also
$h^{1,0}(M)=h^{1,3}(X)=0$ we get from \cite[Cor.\
I:3.7]{illusie83::les+rham+witt} that $k > 0$. By duality
$H^1(M)=\domU_{-k}(-1)$ and if the extension is trivial we wouldn't have that
$M$ is $\cs$-acyclic. Again by acyclicity we get from Lemma \ref{Nk} that $M$ is
isomorphic to $N_1$.
\end{proof}
\end{proposition}
\begin{remark}
Even though we have refrained from proving it, the fact that the de Rham-Witt
cohomology is uniquely determined up to finite ambiguity has an interesting
consequence, namely that contrary to the case of supersingular K3-surfaces we
can not hope to get a classification of supersingular (for some definition of
that term that would include our examples) Calabi-Yau threefolds through
cohomology (at least not the de Rham-Witt cohomology) as we have positive
dimensional families.
\end{remark}
\end{subsection}
\end{section}
\bibliography{preamble,abbrevs,alggeom,algebra,ekedahl}
\bibliographystyle{pretex}
\end{document}